\documentclass[11pt]{article}

\usepackage[margin=1in]{geometry}

\usepackage{amsmath, amssymb, amsfonts, amsthm, mathtools, bm, mathrsfs}

\usepackage{algorithm}
\usepackage{algorithmic}

\usepackage{subcaption}
\usepackage{arydshln}
\usepackage{enumerate}

\usepackage{color}
\usepackage{verbatim}
\usepackage{graphicx}
\graphicspath{{figures/}} 

\usepackage[
    colorlinks=true,
    linkcolor=blue,
    citecolor=blue,
    urlcolor=blue
]{hyperref}
\usepackage[capitalize, nameinlink]{cleveref}
\Crefname{ALC@unique}{Line}{Lines} 

\numberwithin{equation}{section}

\theoremstyle{plain}
\newtheorem{theorem}{Theorem}[section]

\newtheorem{definition}[theorem]{Definition}
\newtheorem{lemma}[theorem]{Lemma}
\newtheorem{proposition}[theorem]{Proposition}
\newtheorem{corollary}[theorem]{Corollary}

\DeclareMathOperator{\Exp}{Exp}
\DeclareMathOperator{\der}{D}
\DeclareMathOperator{\T}{T}
\DeclareMathOperator{\F}{F}
\DeclareMathOperator{\Her}{H}

\DeclareMathOperator{\gl}{GL}
\DeclareMathOperator{\diagm}{diag}
\DeclareMathOperator{\dist}{dist}
\DeclareMathOperator{\vvec}{vec}
\DeclareMathOperator{\nlog}{nlog}
\DeclareMathOperator{\ivg}{ig}

\newcommand{\codecomment}[1]{\hfill {\footnotesize \texttt{// #1}}}

\newcommand{\lblines}[1]{\left\{ \begin{aligned} #1 \end{aligned}\right.}
\newcommand{\realmat}[1]{\mathbb{R}^{#1 \times #1}}

\newcommand{\realset}{\mathbb{R}}
\newcommand{\cmpxmat}[1]{\mathbb{C}^{#1 \times #1}}

\newcommand{\cmpxset}{\mathbb{C}}

\let\origvdots\vdots
\let\origddots\ddots

\newcommand{\svdots}{\scalebox{0.6}{$\origvdots$}}
\newcommand{\sddots}{\scalebox{0.6}{$\origddots$}}

\newcommand{\sbmatrix}[1]{%
  \begingroup
  \let\vdots\svdots
  \let\ddots\sddots
  \left[\begin{smallmatrix}#1\end{smallmatrix}\right]
  \endgroup
}

\newcommand{\im}{\imath}
\newcommand{\zerov}{\mathbf{0}}

\newcommand{\sinc}{\mathop{\mathrm{sinc}}\nolimits}
\newcommand{\cosc}{\mathop{\mathrm{cosc}}\nolimits}
\newcommand{\cotc}{\mathop{\mathrm{cotc}}\nolimits}
\newcommand{\dexp}{\der\exp}
\newcommand{\dexpof}[1]{\dexp\left(#1\right)}
\newcommand{\dexpinvof}[1]{\left(\dexp\left(#1\right)\right)^{-1}}
\newcommand{\so}{\mathcal{SO}}

\newcommand{\skewm}{\mathbf{Skew}}
\newcommand{\skewf}{\mathcal{S}}

\newcommand{\core}{\mathcal{C}_{\Theta}}
\newcommand{\coreop}{\mathcal{C}}
\newcommand{\dskew}{\mathcal{L}}
\newcommand{\conjset}{\fraks}
\newcommand{\poly}[2]{[#1]_{#2}}
\newcommand{\ball}{\mathcal{B}}

\newcommand{\calm}{\mathcal{M}}
\newcommand{\calz}{\mathcal{Z}}
\newcommand{\fraks}{\mathfrak{S}}
\newcommand{\frakz}{\mathfrak{Z}}

\title{The Exponential of Skew-Symmetric Matrices: A Nearby Inverse and Efficient Computation of Derivatives\thanks{
    This work was supported by the National Natural Science Foundation of China (No. 12371311), the Natural Science Foundation of Fujian Province (No. 2023J06004), the Fundamental Research Funds for the Central Universities (No. 20720240151), and the Fonds de la Recherche Scientifique-FNRS under Grant no T.0001.23.
}}

\author{Zhifeng Deng\thanks{School of Mathematical Sciences, Xiamen University, Xiamen, China
  (\url{zhifengdeng@xmu.edu.cn}).}
\and P.-A. Absil\thanks{ICTEAM Institute, UCLouvain, Louvain-la-Neuve, Belgium
  (\url{pa.absil@uclouvain.be}).}
\and Kyle A. Gallivan\thanks{Department of Mathematics, Florida State University, Tallahassee, USA (\url{kgallivan@fsu.edu}).}
\and Wen Huang\thanks{Corresponding author. School of Mathematical Sciences, Xiamen University, Xiamen, China
(\url{wen.huang@xmu.edu.cn}).}
}

\begin{document}
\maketitle

\begin{abstract}
    The matrix exponential restricted to skew-symmetric matrices has numerous applications, notably in view of its interpretation as the Lie group exponential and Riemannian exponential for the special orthogonal group. We characterize the invertibility of the derivative of the skew-restricted exponential, thereby providing a simple expression of the tangent conjugate locus of the orthogonal group. In view of the skew restriction, this characterization differs from the classic result on the invertibility of the derivative of the exponential of real matrices. Based on this characterization, for every skew-symmetric matrix $A$ outside the (zero-measure) tangent conjugate locus, we explicitly construct the domain and image of a smooth inverse---which we term \emph{nearby logarithm}---of the skew-restricted exponential around $A$. This nearby logarithm reduces to the classic principal logarithm of special orthogonal matrices when $A=\zerov$. The symbolic formulae for the differentiation and its inverse are derived and implemented efficiently. The extensive numerical experiments show that the proposed formulae are up to $3.9$ times and $3.6$ times faster than the current state-of-the-art robust formulae for the differentiation and its inversion, respectively. 
\end{abstract}



\section{Introduction}\label{sec:introduction}

The matrix exponential arises in many important applications, e.g., quantum theory~\cite{wilcox1967exponential}, economics and statistics~\cite{chan2003generalized, chen2001analytic}, computer vision~\cite{turaga2011statistical}, and neural network training~\cite{zhou2019continuity}. Much research effort has been put into the map $\exp:\cmpxmat{n}\to \gl_n(\cmpxset)$, where $\gl_n(\cmpxset)$ and $\gl_n(\realset)$ are the general linear groups consisting of invertible complex matrices and real matrices, respectively. Al-Mohy and Higham~\cite{al2009exponential,al2010new} investigated the exponential of arbitrary complex matrices, and Najfeld and Havel~\cite{najfeld1995} also investigated the exponential of diagonalizable matrices.

The restriction of the matrix exponential to the skew-symmetric matrices deserves special consideration. From here onward, unless otherwise stated, $\exp$ stands for
\begin{equation*}
    \begin{aligned}
        \exp:\skewm_n\to \so_n& \text{, where } \lblines{
            &\skewm_n := \{A\in\realmat{n}:A+A^{\T}=\zerov\},\\
            &\so_n:=\{Q\in\realmat{n}:Q^{\T}Q =I_n, \det(Q) = 1\}.\\
        }
    \end{aligned}
\end{equation*}
As a Lie group, $\so_n$ arises in various applications, e.g., the orthogonal weights in neural network training~\cite{wang2020orthogonal} and image representation~\cite{qi2021orthogonalmoment}. Meanwhile, $\skewm_n$ also arises in algorithms on $\so_n$. For example, the continuous orthonormalization algorithm~\cite{dieci2000continuous, meyer1986continuous} applied to various PDE problems~\cite{barker2018evans, borzsak1996lyapunov, dieci1997compuation} converts $\dot{Y}(t) = C(t)Y(t)$ into a new PDE of the form $\dot{Q}(t) = Q(t)A(t)$, where $Q(t)$ denotes the Q factor of the QR decomposition $Y(t) = Q(t)R(t)$ and $A(t)$ is constrained to $\skewm_n$. Furthermore, developments in orthogonal constrained manifolds have led to equations of the form $\exp(A) = Q$, e.g., the geodesics~\cite{edelman1998} and extremal curves~\cite{jurdjevic2020extremal} on Stiefel and Grassmann manifolds are characterized by such an equation with partially known $A\in\skewm_n$ and $Q\in\so_n$. Solving these PDE problems and constrained equations on manifolds not only demands significant computations involving the matrix exponential and/or its differentiation, but also requires careful consideration of the invertibility of the matrix exponential and its differentiation. In particular, $\log(\exp(A))$ fails to return $A$ whenever $A \notin \mathrm{image}(\log) = \{A \in\skewm_n :\|A\|_2 < \pi\}$,  where $\log$ denotes the principal logarithm on $\so_n$. As for the differentiation, it is known that, when $n\geq 3$, the derivative of $\exp$ is not everywhere an invertible linear map; otherwise $\exp$ would be an open map, contradicting~\cite[Proposition~8.3]{absil2025ultimate}. These limitations explain why, e.g., the algorithm for computing endpoint Stiefel geodesics~\cite[Algorithm~4]{zimmermann2021computing} assumes sufficiently close points.

However, there is less literature that exploits the skew-symmetry in the exponential map to meet the growing needs specific to $\so_n$. From the computational point of view, there are Rodrigues-like formulae proposed in~\cite{gallier2003computing} and Pad{\'e} approximant proposed in~\cite{cardoso2010exponentials} that accelerate $\exp:\skewm_n\to\so_n$ and the principal logarithm.

This paper presents novel and efficient Schur-based formulae in real arithmetic for $\dexpof{A}:T_A\skewm_n\to T_Q\so_n$ and $\dexpinvof{A}:T_Q\so_n\to T_A\skewm_n$ where $A\in\skewm_n,Q=\exp(A)\in\so_n$. As an important theoretical implication, this paper proposes a notion of nearby matrix logarithm on $\so_n$ as a local smooth inversion to $\exp:\skewm_n\to \so_n$. This nearby matrix logarithm subsumes the principal logarithm, and it naturally works beyond the image of the principal logarithm.

\subsection{Related Work}

The idea of utilizing the differentiation $\dexpof{A}$ and its inversion for computing the logarithm of nearby matrices on the general linear group $\gl_n$ was first proposed in~\cite{dieci1999real}, in which Dieci et al. proposed a quadrature-based formula for $\dexpof{A}:\realmat{n}\to\realmat{n}$ for $A$ in the image of the principal logarithm. Currently, the standard approach of computing the matrix exponential and the associated objects uses the order-$[r/s]$ Pad{\'e} approximant $f(X)\approx \poly{r/s}{f}(X) := (h_s(X))^{-1} (g_r(X))$,
where the polynomials $g_r(X) = \sum_{i=0}^r \alpha_iX^i$ and $h_s(X)=\sum_{i=0}^s \beta_iX^i$ are determined by $f$. Additionally, a scaling process $f(X)\approx (f(X / 2^{\sigma}))^{2^{\sigma}}$ with an integer $\sigma > 0$ is necessary for $X$ with significant matrix norm. In particular, the scaling and squaring method computes the matrix exponential $\poly{r/r}{\exp}$ with $r = 3,\ldots,13$, see~\cite{higham2005scaling}, which requires up to $6 + \sigma$ matrix multiplications, one LU decomposition and two calls to the LU-solver with different $n\times n$ matrices.

Al-Mohy and Higham~\cite{al2009exponential} further derived $\dexpof{A}:\cmpxmat{n}\to \cmpxmat{n}$ based on $\poly{r/r}{\exp}$, which requires up to $(19 + 3  \sigma)$ matrix multiplications, one LU decomposition and two calls to the LU-solver with different $n\times n$ matrices. Later, Al-Mohy et al.~\cite{al2013logarithm} derived the formula that differentiates the principal logarithm $\der\log:\cmpxmat{n}\to \cmpxmat{n}$ using $\poly{r/r}{\log}, r = 1,\ldots,7$, which requires a Schur decomposition of $Q$, six matrix multiplications, $r$ calls to a Sylvester equation solver, $\sigma$ calls to LU decompositions and $2\sigma$ calls to the LU-solver with different $n\times n$ matrices. Note that $\der\log(Q)$ only agrees with $\dexpinvof{A}$ when $A$ is in the image of the principal logarithm. These formulae adopt real arithmetic with real matrices. However, the skew-symmetry not only allows more aggressive computation techniques, but also implies unexploited structures in the Pad{\'e} approximant for general matrices.

Najfeld and Havel~\cite{najfeld1995} proposed a formula for $\dexpof{A}:\cmpxmat{n}\to\cmpxmat{n}$ at a diagonalizable $A = V\Lambda V^{-1}$, which requires an eigendecomposition and four matrix multiplications. For a skew-symmetric $A$ that is diagonalizable with purely imaginary eigenvalues and the non-real eigenvectors, the formula only works in complex arithmetic. A contribution of this paper is to further develop the complex formula in the instance of a skew-symmetric matrix such that the computation is performed in real arithmetic and the geometric implications are exploited.

The local invertibility in differentiating $\exp:\cmpxmat{n}\to \gl_n(\cmpxset)$ or $\exp:\realmat{n}\to \gl_n(\realset)$ is also investigated from the Lie group perspective in~\cite[Sec.~1.2]{rossmann2006lie}, which is established at $A$ if there is no pair of eigenvalues that differ by a multiple of $2\pi\im$. Dieci et al.~\cite{dieci1999real} consider the problem of obtaining the logarithm of $A + E$ once the logarithm of $A$ is known and $E$ is a matrix of small norm. While these results are applicable to the $\exp:\skewm_n\to\so_n$, finer and more restrictive geometric structures with stronger conclusions on both $\skewm_n$ and $\so_n$ have not been thoroughly explored in the literature.

Recently, Mataigne and Gallivan~\cite{mataigne2024eigenvalue} presented an accelerated formula for the Schur decomposition of skew-symmetric matrices, which is a fundamental primitive that guarantees the efficiency of the work presented in this paper. Gallier and Xu~\cite{gallier2003computing} also proposed Rodrigues-like formulae for accelerating $\exp:\skewm_n\to \so_n$ and the principal logarithm $\log:\so_n\to\skewm_n$. Unlike the Pad{\'e} approximant that can be reused in $\dexpof{A}$ and $\der\log(Q)$, or the Schur decomposition of $A$ that can be reused in the proposed new formulae, the Rodrigues-like formulae in~\cite{gallier2003computing} do not compute objects necessary for differentiation or its inverse, making them less applicable to the problem of interest in this paper.

\subsection{Contributions}

Let $\dskew_A:\skewm_n\to\skewm_n$ be defined by:
\begin{equation}\label{eq:dexp-skew}
    \begin{aligned}
        \dexpof{A}: T_A\skewm_n&\to T_Q\so_n: &X &\mapsto \dexpof{A}[X]= Q \dskew_A(X).
    \end{aligned}
\end{equation}
The main contribution of this paper are the following:
\begin{enumerate}
    \item The invertibility of $\dexpof{A}$ on skew-symmetric matrices is characterized and discussed via $\dskew_A$.
    \item Fast symbolic formulae for $\dskew_A$ and $\dskew_A^{-1}$, which characterize $\dexpof{A}$ and $\dexpinvof{A}$, are proposed and implemented (in \texttt{C++}). In contrast with~\cite{najfeld1995}, the computation of $\dexpof{A}$ is performed in real arithmetic.
    \item An explicit characterization of the nearby matrix logarithm as a local inversion of $\exp:\skewm_n\to \so_n$ is presented, which subsumes the principal logarithm.
\end{enumerate}

\subsection{Organization}

The rest of this paper is organized as follows. \Cref{sec:prelimin-notation} reviews some important statements and introduces the notation used. The symbolic formulae of differentiation and its inversion are presented in \Cref{sec:differentiation}. The algorithmic discussion and complexity analysis are also given there. \Cref{sec:nearby-logarithm} applies the invertibility condition to construct the nearby matrix logarithm on $\so_n$. The numerical experiments on differentiations and their results are presented in \Cref{sec:experiments}.

\section{Preliminaries and Notation}\label{sec:prelimin-notation}

This section introduces the notation on block partitions and some trigonometric quantities associated with the $\sinc$ function widely used in signal processing~\cite{brigham1988fast}. Some preliminaries about the Schur decomposition and eigendecomposition of skew-symmetric matrices are also reviewed. More details of matrix decompositions can be found in Golub and Van Loan~\cite{golub2013matrix}.

\subsection{Even/Odd Dimension and Size-2 Partitions}

The nonzero eigenvalues of $A\in\skewm_n$ form conjugate pairs $\pm\theta\im$ for some $\theta\in\realset$. When $n=2m$, the eigenvalues of $A$ are of the form $\theta_1\im, -\theta_1\im, \dots, \theta_m\im, -\theta_m\im$ with $\Theta:=(\theta_1,\dots,\theta_m) \in \realset^m$. When $n = 2m+1$, there is an extra zero eigenvalue. For the even case, $m = n/2$ is sufficient for all derivations but when $n=2m+1$ is odd, the integer $k = m+1 = (n+1)/2$ also appears frequently. Therefore, the following notation of $n, m$ and $k$ are assumed throughout this paper.

The integer $n$ refers to the dimension of $A\in\skewm_n$. The integers $m$ and $k$ are determined by taking the floor and ceiling of $n/2$
\begin{equation*}
    \begin{aligned}
        &m := \left\lfloor \frac{n}{2} \right\rfloor = \begin{cases}
            n / 2 & n \text{ is even},\\
            (n - 1) / 2 & n \text{ is odd},\\
        \end{cases}
        &k := \left\lceil \frac{n}{2} \right\rceil = \begin{cases}
            n / 2 & n \text{ is even},\\
            (n + 1) / 2 & n \text{ is odd}.\\
        \end{cases}
    \end{aligned}
\end{equation*}

By grouping conjugate pairs of eigenvalues in $A\in\skewm_n$ appropriately, both the exponential $\exp(A)$ and the differentiation $\dexpof{A}[X]$ can be expressed in terms of blocks of size two, which are denoted by the squared brackets as follows.

For a squared $n\times n$ matrix $X$, the subscript with one bracketed index $X_{[i]}, i \leq m$ denotes the $i$th pair of columns in $X$. When $n = 2m+1$ is odd, $X_{[k]}, k = m+1$ denotes the left-over column. Similarly, the subscript with two bracketed indices $X_{[i,j]}, i,j\leq m$ denotes the $i,j$th $2\times 2$ block in $X$ while the $X_{[i,k]}, X_{[k,j]}$ and $X_{[k,k]}$ with $i,j\leq m, k = m+1$ denote the blocks in the left-over row, column, and diagonal as demonstrated in the $n=5$ case.

When $n = 5, m= 2, k =3$, the block partitions are given by
\begin{equation*}
    \begin{aligned}
        \left[\begin{array}{cc:cc:c}
            X_{11} & X_{12} & X_{13} & X_{14} & X_{15}\\
            X_{21} & X_{22} & X_{23} & X_{24} & X_{25}\\
            \hline
            X_{31} & X_{32} & X_{33} & X_{34} & X_{35}\\
            X_{41} & X_{42} & X_{43} & X_{44} & X_{45}\\
            \hline
            X_{51} & X_{52} & X_{53} & X_{54} & X_{55}\\
        \end{array}\right] \begin{aligned}
            &=\left[\begin{array}{c:c:c}
                X_{[1,1]} & X_{[1, 2]} & X_{[1, 3]}\\
                \hline
                X_{[2,1]} & X_{[2, 2]} & X_{[2, 3]}\\
                \hline
                X_{[3,1]} & X_{[3, 2]} & X_{[3, 3]}\\
            \end{array}\right]\text{ in blocks,}\\
            &\text{ or }=\left[\begin{array}{c:c:c}
                X_{[1]} & X_{[2]} & X_{[3]}\\
            \end{array}\right] \text{ in columns.}
        \end{aligned}\\
    \end{aligned}
\end{equation*}

\subsection{Schur Decomposition}

The (real) Schur decomposition of $X\in\realmat{n}$ transforms the matrix into a block upper triangular $T$ with the orthogonal Schur vectors $R$ as $X = RTR^{\T}$, where the diagonal blocks in $T$ are $2\times 2$ or $1\times 1$. When $X$ is a normal matrix (i.e., $X^{\T}X = XX^{\T}$), the block upper triangular $T$ becomes block diagonal where the $2\times 2$ blocks share the form of $\sbmatrix{a & -b\\b & a},a,b\in\realset$.

\begin{lemma}\label{lemma:schur-decomposition-skew-so}\cite[Corollary 2.5.11]{horn2012matrix}
        There exists a Schur decomposition for any $A = R\Xi R^{\T}\in\skewm_n$ or $Q = RER^{\T}\in\so_n$, such that, for all $i\leq m$, $\Xi_{[i,i]} = \sbmatrix{0 & -\theta_i\\\theta_i & 0}$ with $\theta_i\in\realset$, or $E_{[i,i]} = \sbmatrix{c_i & -s_i\\s_i & c_i}$ with $c_i^2+s_i^2 = 1$. For $n = 2m+1$, the left-over $1\times 1$ diagonal is $E_{[k,k]} = 1$ or $\Xi_{[k,k]} = 0$.
\end{lemma}
\begin{proof}
    When $A=R\Xi R^{\T}\in\skewm_n$, $\Xi$ is skew-symmetric, so are the diagonal blocks $\Xi_{[i,i]},i\leq k$. All $2\times 2$ skew-symmetric matrices share the form, $\sbmatrix{ 0 & -\theta_i\\\theta_i & 0}$, and the zero matrix $0$ is the only $1\times 1$ skew-symmetric matrix. Similarly, the form of the diagonal blocks $E_{[i,i]}, i\leq k$ is given by the special orthogonality.
\end{proof}

The Schur decompositions of $A = R\Xi R^{\T}\in\skewm_n$ and $Q = RER^{\T}\in\so_n$ in the forms given by \cref{lemma:schur-decomposition-skew-so} are denoted as $(\Xi, R)$ and $(E, R)$ respectively, where $\Xi$ is reserved for skew-symmetric matrices and $E$ is reserved for special orthogonal matrices. Moreover, $\Theta := (\theta_1,\ldots,\theta_m)\in\realset^m$ determined by $\Xi$ is referred to as the set of angles of $A$ (under $R$), and the redundant tuple $(\Theta, \Xi, R)$ also denotes $A=R\Xi R^{\T}\in\skewm_n$.

The notion of angles of $A\in\skewm_n$ enables important statements such as \cref{thm:invertibility-skew-symm} and \cref{prop:distance-to-subset}. Although a normal matrix $X$ has multiple Schur decompositions, there exists an invariance group of $T = R^{\T}XR$ defined as
\begin{equation}\label{eq:invariance-group}
    \ivg(T) := \{P\in\realmat{n}:P^{\T}P = I_n, P^{\T}TP = T\}.
\end{equation}
If $(\Xi,R)$ is a Schur decomposition of $A \in \skewm_n$, then $(\Xi,\tilde{R})$ is also a Schur decomposition if and only if $\tilde{R} = RP$ with $P\in\ivg(\Xi )$. Note that the invariance group preserves the order and multiplicities of the angles $\Theta\in\realset^m$ determined in $\Xi$. The invariance group \cref{eq:invariance-group} is addressed and discussed in~\cite{mataigne2024eigenvalue}, and a property of invariance groups given in \cref{lemma:invariance-group} is important to the development in \Cref{sec:nearby-logarithm}.

\begin{lemma}\label{lemma:invariance-group}
    Consider two block diagonal normal matrices $Y=\diagm(Y_1,\ldots, Y_l)$ and $Z=\diagm(Z_1,\ldots, Z_l)$ where $Y_i$ and $Z_i$ are normal blocks with the same dimension of $2\times 2$ or $1\times 1$, then $\ivg(Y) = \ivg(Z)$ if $Y_i = Y_j \Leftrightarrow Z_i = Z_j,\forall i,j$.
\end{lemma}
\begin{proof}
    According to~\cite[Lemma A.1$\sim$A.3]{mataigne2024eigenvalue}, the invariance group $\ivg(Y)$ of the normal block diagonal $Y$ only depends on the multiplicities of the eigenvalues of $Y$. Since the condition $Y_i = Y_j \Leftrightarrow Z_i = Z_j$ preserves such multiplicities of $Y$ in the eigenvalues of $Z$, it holds that $\ivg(Y) = \ivg(Z)$.
\end{proof}

For example, let $G$ be a $2\times 2$ rotation, i.e., $G = \sbmatrix{c & -s\\s & c}$, $c^2+s^2 = 1$. It holds that $G^{\T}\sbmatrix{0 & -x\\x & 0}G = \sbmatrix{0 & -x\\x & 0}$ for all $x\in\realset$. Hence $R_{[i]}G$ consists of Schur vectors associated with the same angle, as $R_{[i]}\Xi_{[i,i]}R_{[i]}^{\T} = (R_{[i]}G)\Xi_{[i,i]}(R_{[i]}G)^{\T}$. On the other hand, reflecting $R_{[i]}$ with any $2\times 2$ reflection $H = \sbmatrix{s & c\\ c& -s}$ yields a sign-flip to $\theta_i$, as $H^{\T}\sbmatrix{0 & -x\\x & 0}H = \sbmatrix{0 & x\\-x & 0}$ for all $x\in\realset$. While the angles can always be converted to nonnegative values, allowing for negative angles is preferable as it makes it possible to smoothly assign the angles $\Theta(t)$ in $\realset^m$ of a smoothly varying $A(t)$ in $\skewm_n$, as shown in \cref{fig:st-curve}. Observe that the conditions in \cref{thm:invertibility-skew-symm} and \cref{prop:distance-to-subset} are invariant by any sign-flip of angles.

\subsection{Eigendecomposition}

An eigendecomposition of a skew-symmetric $A$ is constructed by a Schur decomposition $(\Theta, \Xi, R)$ as given in \cref{prop:spectral-decomposition}.

\begin{proposition}\label{prop:spectral-decomposition}
    Let $(\Theta, \Xi, R)$ be a real Schur decomposition of $A\in\skewm_n$, $U := \frac{\sqrt{2}}{2} \sbmatrix{
        1 & 1\\
        -\im&\im
    }$, and $U_n$ be the following:
    \begin{equation}\label{eq:block-diagonal-Un}
        U_n:= \begin{cases}
            \diagm\left(U,\ldots, U\right) & n = 2m,\\
            \diagm\left(U,\ldots, U, 1\right) & n = 2m+1.\\
        \end{cases}
    \end{equation}
    Then, $(\Theta, \Xi, R)$ yields an eigendecomposition $(\Lambda, V=RU_n)$ of $A$ as follows
    \begin{equation}\label{eq:schur-to-spectral}
        \begin{aligned}
            &\begin{cases}
                \left(\begin{aligned}
                    &\Lambda = \diagm\left(-\theta_1\im , \theta_1\im, \ldots, -\theta_m\im, \theta_m\im \right),\\
                    &V=RU_n = \begin{bmatrix}
                        R_{[1]}U & \cdots & R_{[m]}U
                    \end{bmatrix}\\
                \end{aligned}\right)& n = 2m,\\
                \left(\begin{aligned}
                    &\Lambda = \diagm\left(-\theta_1\im , \theta_1\im, \ldots, -\theta_m\im, \theta_m\im , 0\right),\\
                    &V=RU_n = \begin{bmatrix}
                        R_{[1]}U & \cdots & R_{[m]}U & R_{[m+1]}
                    \end{bmatrix}\\
                \end{aligned}\right)& n = 2m+1.\\
            \end{cases}\\
        \end{aligned}
    \end{equation}
\end{proposition}
\begin{proof}
    This follows from the identity $U\sbmatrix{
        -x\im & 0\\
        0 & x\im
    }U^{\Her} = \sbmatrix{
        0 & -x\\
        x & 0
    },\forall x\in\realset$.
\end{proof}

\subsection{Trigonometric Notation}

There are three trigonometric quantities that frequently appear in the discussion. These are given in \cref{def:trig-func}.
\begin{definition}\label{def:trig-func}
    The $\sinc$, $\cosc$ and $\cotc$ functions are given by
    \begin{equation}\label{eq:trigonometries}
        \begin{aligned}
            \sinc(x) &= \begin{cases}
                1 & x = 0,\\
                \sin(x) \slash x & x\neq 0,\\
            \end{cases}&
            \cosc(x) = \begin{cases}
                0 & x = 0,\\
                \left(\cos(x)-1\right) \slash x & x \neq 0,\\
            \end{cases}\\
            \cotc(x) &= \begin{cases}
                1 & x = 0,\\
                \text{DNE} & x = \pm \pi,\ldots,\\
                x  \cot(x) & \text{otherwise}.\\
            \end{cases}& (\text{DNE stands for ``does not exist''.})
        \end{aligned}
    \end{equation}
\end{definition}

\section{Differentiation of the Matrix Exponential and its Inversion}\label{sec:differentiation}

This section describes the derivation of the formulae, in real arithmetic, that compute the differentiation of the matrix exponential, \cref{eq:dskew-skew}, and its inversion, \cref{eq:dexp-skew-inv}. The complexity analyses, \cref{tab:dexp-ops} and \cref{tab:dlog-ops}, imply significantly improved computational performance of the proposed new formula, which is further supported by the experiments in \Cref{sec:experiments}. From the theoretical point of view, the formula \cref{eq:dexp-skew-inv} leads to the invertibility condition \cref{eq:invertible-condition} specific for $\dexpof{A}:\skewm_n\to T_Q\so_n$.

Recall that the differentiation of the matrix exponential at $A\in\skewm_n$ with $Q = \exp(A)\in\so_n$ in~\eqref{eq:dexp-skew}, and its inversion, act between the tangent spaces
\begin{equation*}
    \begin{aligned}
        \lblines{
        &\dexpof{A}:T_A\skewm_n\to T_Q\so_n,\\
        &\dexpinvof{A}:T_Q\so_n \to T_A\skewm_n,
    } & \text{ where } &\lblines{
        &T_A\skewm_n = \skewm_n,\\
        &T_Q\so_n = \{Q X:X\in\skewm_n\}.
    }
    \end{aligned}
\end{equation*}
The similar structure of $T_A\skewm_n$ and $T_Q\so_n$ leads to a linear transformation of $\skewm_n$ determined by $A$, that was introduced in~\cref{eq:dexp-skew}, in the form of
\begin{equation}\label{eq:linear-skew}
    \begin{aligned}
        \dskew_A:\skewm_n &\to \skewm_n:&X &\mapsto Q^{\T} \dexpof{A}[X] = \dskew_A(X),
    \end{aligned}
\end{equation}
such that both the differentiation and its inverse are expressed as
\begin{equation*}
    \lblines{
        &\dexpof{A}[X] = Q \dskew_A(X),\\
        &\dexpinvof{A}[Q X] = \dskew_A^{-1}(X),
    } \forall X\in\skewm_n.
\end{equation*}

As linear transformations of $\skewm_n\subset \realmat{n}$, $\dskew_A$ and $\dskew_A^{-1}$ must admit formulae that perform all computations in real arithmetic.

\subsection{Linear Transformation in Complex Arithmetic}

For all diagonalizable $A$, there exists a Daletski$\breve{\mathrm{\i}}$--Kre$\breve{\mathrm{\i}}$n \footnote{Also spelled Daleckii--Krein in the literature.} formula~\cite{daletskii1965integration, najfeld1995} that symbolically computes $\dexpof{A}:\cmpxmat{n}\to\cmpxmat{n}$ based on an eigendecomposition. This complex formula is ideal for further developing the real formula of $\dskew_A$ as (i) $A\in\skewm_n$ is diagonalizable, (ii) an eigendecomposition of $A$ can be obtained from a real Schur decomposition of $A$ as in~\eqref{eq:schur-to-spectral}, and (iii) symbolic formulae help simplify the discussion of invertibility.


\begin{lemma}[Daletski$\breve{\mathrm{\i}}$--Kre$\breve{\mathrm{\i}}$n Formula~\cite{daletskii1965integration, najfeld1995}]\label{lemma:dexp-diag}
    The differentiation $\dexpof{A}$ at a diagonalizable $A\in \cmpxmat{n}$ along a perturbation $X\in \cmpxmat{n}$ is given by
    \begin{equation}\label{eq:dexp-diag}
        \dexpof{A}[X] = V\left((V^{-1}X V)\odot \Psi\right)V^{-1},
    \end{equation}
    where $(\Lambda,V)$ is an eigendecomposition of $A$, $\odot$ is the Hadamard product that performs entry-wise multiplication, and the matrix $\Psi\in\cmpxmat{n}$ has the entries
    \begin{equation}\label{eq:dexp-diag-para}
        \psi_{ij} = \begin{cases}
                \left(e^{\lambda_j} - e^{\lambda_i}\right)\big/\left(\lambda_j - \lambda_i\right) & \lambda_i\neq \lambda_j,\\
                e^{\lambda_i} & \lambda_i = \lambda_j.\\
            \end{cases}
    \end{equation}
\end{lemma}

\begin{proposition}\label{prop:extract-linear-skew}
    The operator $\dskew_A$ in $\dexpof{A}[X] =: \exp(A) \dskew_A(X)$ is
    \begin{equation}\label{eq:linear-diag}
        \dskew_A(X) = V\left((V^{-1}X V)\odot \Phi\right)V^{-1},
    \end{equation}
    where the matrix $\Phi = \exp(-\Lambda)\Psi$ has the entries
    \begin{equation}\label{eq:linear-diag-para}
        \phi_{ij} = \begin{cases}
                \left(e^{\lambda_j-\lambda_i} - 1\right)\big/ \left(\lambda_j - \lambda_i\right) & \lambda_i\neq \lambda_j,\\
                1 & \lambda_i = \lambda_j.\\
            \end{cases}
    \end{equation}
\end{proposition}
\begin{proof}
    For a diagonalizable $A = V\Lambda V^{-1}$, there is
    \begin{equation*}
        \begin{aligned}
            \dskew_A(X) &= \exp(-A) \dexpof{A}[X] = \left(V\exp(-\Lambda)V^{-1}\right) \left(V\left((V^{-1}X V)\odot \Psi\right)V^{-1}\right)\\
            &= V\left(\exp(-\Lambda) \left((V^{-1}X V)\odot \Psi\right)\right)V^{-1} = V\left((V^{-1}X V)\odot (\exp(-\Lambda)\Psi)\right)V^{-1}.
        \end{aligned}
    \end{equation*}
    The last equality follows from the observation that the $(i,j)$th entry of $N = \Sigma (M\odot H)$ with arbitrary matrices $M,H$ and a diagonal $\Sigma = \diagm(\sigma_1,\ldots, \sigma_n)$ is given by
    \begin{equation*}
        N_{ij} = \sigma_i M_{ij} H_{ij} = M_{ij}  \left(\sigma_i H_{ij}\right)
    \end{equation*}
    where $\sigma_i H_{ij}$ is the $(i,j)$th entry of $\Sigma  H$, s.t. $\Sigma (M\odot H) = M\odot(\Sigma  H)$.
\end{proof}

\subsection{Core Map in Real Arithmetic}\label{subsec:core-map}

As mentioned in \cref{prop:spectral-decomposition}, the eigenvectors $V$ can be expressed as the product of the Schur vectors $R$ and the unitary matrix $U_n$. Substituting $V = RU_n$ into~\eqref{eq:dexp-diag} and grouping the actions of $R$ together (in real arithmetic) yields a core map in complex arithmetic that excludes all influences from $R$, which is
\begin{equation*}
    \begin{aligned}
        &\dskew_A(X) = V\left((V^{-1}X V)\odot \Phi\right)V^{-1} = R\core(R^{\T}XR)R^{\T}
    \end{aligned}
\end{equation*}
where $\core:\skewm_n\to\skewm_n$ denotes the \emph{core map} as follow
\begin{equation}\label{eq:core-map}
    \begin{aligned}
        \core:\skewm_n&\to\skewm_n:&M &\mapsto N := U_n\left((U_n^{\Her} M U_n) \odot \Phi\right)U_n^{\Her}.
    \end{aligned}
\end{equation}
In view of \cref{eq:linear-skew}, the differentiation is given by $\dexpof{A}[X] = Q  \left(R \core(R^{\T}  X   R)  R^{\T}\right)$ for all $A, X\in\skewm_n$ where the computations other than the core map are multiplications of real matrices. It remains to derive a real formula for the core map.

For all $M\in\skewm_n$, applying the core map to get $N = \core(M)$ consists of the following three sequential operations
\begin{equation*}
    \begin{matrix}
        &M \mapsto M' := U_n^{\Her}M U_n, & M'\mapsto M'' :=M'\odot \Phi, & M'' \mapsto N:= U_n M'' U_n^{\Her}.
    \end{matrix}
\end{equation*}
For the first and the last congruences with $U_n$ and $U_n^{\Her}$, the block diagonal $U_n$ in~\eqref{eq:block-diagonal-Un} yields block multiplications in $M_{[i,j]}$ and $M''_{[i,j]},i, j\leq k$, e.g., $M'$ with $n = 5$ is
\begin{equation*}
    M' =\sbmatrix{
        U^{\Her}M_{[1,1]} U & U^{\Her}M_{[1,2]}U & U^{\Her}M_{[1,3]}\\
        U^{\Her}M_{[2,1]} U & U^{\Her}M_{[2,2]}U & U^{\Her}M_{[2,3]}\\
        M_{[3,1]} U & M_{[3,2]}U & M_{[1,3]}\\
    },
\end{equation*}
and similarly $N = \diagm(U, U, 1)  M'' \diagm(U^{\Her}, U^{\Her}, 1)$.

The Hadamard product in $M''= M'\odot \Phi$ can be grouped as block Hadamard products $M''_{[i,j]}= M'_{[i,j]}\odot \Phi_{[i,j]},\forall i,j\leq k$. In view of \cref{eq:linear-diag-para} and \cref{eq:schur-to-spectral} the block $\Phi_{[i,j]}, i,j\leq k$ is fully determined by $\theta_i$ and/or $\theta_j$. For $i,j\leq m$, $\Phi_{[i,j]}$ is $2\times 2$ and it is determined by the eigenvalues $\pm\theta_i\im,\pm\theta_j\im$. When $\theta_i\neq \theta_j$, it holds that
\begin{equation}\label{eq:a2d-formulae}
    \begin{aligned}
        \Phi_{[i,j]} &= \begin{bmatrix}
            \frac{e^{\lambda_{2j-1}-\lambda_{2i-1}}-1}{\lambda_{2j-1}-\lambda_{2i-1}} &\frac{e^{\lambda_{2j}-\lambda_{2i-1}}-1}{\lambda_{2j}-\lambda_{2i-1}}\\
            \frac{e^{\lambda_{2j-1}-\lambda_{2i}}-1}{\lambda_{2j-1}-\lambda_{2i}} &\frac{e^{\lambda_{2j}-\lambda_{2i}}-1}{\lambda_{2j}-\lambda_{2i}}
        \end{bmatrix}= \begin{bmatrix}
            \frac{e^{(-\theta_j+\theta_i)\im}-1}{(-\theta_j+\theta_i)\im} & \frac{e^{(\theta_j+\theta_i)\im}-1}{(\theta_j+\theta_i)\im}\\
            \frac{e^{(-\theta_j-\theta_i)\im}-1}{(-\theta_j-\theta_i)\im} & \frac{e^{(\theta_j-\theta_i)\im}-1}{(\theta_j-\theta_i)\im}\\
        \end{bmatrix}\\
        &= \begin{bmatrix}
            \frac{\cos(-\theta_j+\theta_i) - 1 + \left(\sin(-\theta_j+\theta_i)\right)\im}{(-\theta_j+\theta_i)\im}& \frac{\cos(\theta_j+\theta_i) - 1 + \left(\sin(\theta_j+\theta_i)\right)\im}{(\theta_j+\theta_i)\im}\\
            \frac{\cos(-\theta_j-\theta_i) - 1 + \left(\sin(-\theta_j-\theta_i)\right)\im}{(-\theta_j-\theta_i)\im}& \frac{\cos(\theta_j-\theta_i) - 1 + \left(\sin(\theta_j-\theta_i)\right)\im}{(\theta_j-\theta_i)\im}\\
        \end{bmatrix}\\
        &:=\begin{bmatrix}
            a-b\im & c-d\im\\
            c+d\im & a+b\im
        \end{bmatrix} \text{ where }\lblines{
            a &= \sinc(\theta_i-\theta_j), &b = \cosc(\theta_i-\theta_j)\\
            c &= \sinc(\theta_i+\theta_j), &d = \cosc(\theta_i+\theta_j)\\
        }.\\
    \end{aligned}
\end{equation}
When $\theta_i=\theta_j$, the diagonals in $\Phi_{[i,j]}$ are replaced by one, which is subsumed in the above $a, \ldots, d$ formulae with $a = \sinc(0) = 1$ and $b = \cosc(0) = 0$.

When $n=2m+1$ is odd, the $2\times 1$ block $\Phi_{[j,k]}$ and the $1\times 2$ block $\Phi_{[k,j]}$, $j\leq m, k = m+1$ are determined by the eigenvalues $\pm\theta_j\im$ and $0$ in the form of
\begin{equation}\label{eq:rs-formaule}
    -\Phi_{[j, k]} = \Phi_{[k,j]}^{\T} = \sbmatrix{-(e^{-\theta_j\im}-1)\big/\theta_j\im \\ (e^{\theta_j\im}-1)\big/\theta_j\im} := \begin{bmatrix}
        r+s\im \\ r-s\im
    \end{bmatrix} \text{ where }\lblines{
        r &= \sinc(\theta_j),\\
        s &= \cosc(\theta_j).\\
    }
\end{equation}
Combining these computations in blocks (that can be computed in-place), the core map $\core$ is written as a collection of linear maps $M_{[i,j]}\mapsto M'_{[i,j]}\mapsto M''_{[i,j]}\mapsto N_{[i,j]}$ for all $i,j\leq k$. The real matrix expressions of the linear maps $M_{[i,j]}\mapsto N_{[i,j]}$ for $i,j\leq m$ are given in \cref{prop:core-action}, which share the forms of $\mathfrak{X}(a,b,c,d)\in\realmat{4}$ and the matrix $\mathfrak{Y}(r,s)\in\realmat{2}$ for real numbers $a,b,c,d,r,s\in\realset$ as follows.
\begin{equation}\label{eq:dexp-skew-mat}
    \lblines{
        &\mathfrak{X}(a,b,c,d) := \frac{1}{2}\sbmatrix{
            a+c & -b-d & b - d & a-c\\
            b+d & a+c & -a + c & b-d\\
            -b+d & -a+c & a + c & -b-d\\
            a-c & -b+d & b + d & a+c\\
        },\\
        &\mathfrak{Y}(r,s) := \sbmatrix{
            r & s\\
            -s & r
        }.\\
    }
\end{equation}

\begin{proposition}\label{prop:core-action}
    For all $M\in\skewm_n$, it holds that
    \begin{equation}\label{eq:core-action}
        \sbmatrix{
            N_{[1,1]} & \cdots & N_{[1,k]}\\
            \vdots & \ddots & \vdots\\
            N_{[k,1]} & \cdots & N_{[k,k]}
        } =  N= \core(M) = \sbmatrix{
            \coreop_{[1,1]}(M_{[1,1]}) & \cdots & \coreop_{[1, k]}(M_{[1,k]})\\
            \vdots & \ddots & \vdots\\
            \coreop_{[k,1]}(M_{[k,1]}) & \cdots & \coreop_{[k, k]}(M_{[k,k]})\\
        },
    \end{equation}
    where the linear maps $\coreop_{[i,j]}, i,j\leq k$, admit the following expression. For the vectorized $\vvec(M_{[i,j]})\in\realset^4, i,j\leq m$ and $\vvec(M_{[m+1,j]})\in\realset^2, j\leq m$ when $n = 2m+1$:
    \begin{equation}\label{eq:dskew-skew}
        \begin{aligned}
            &\begin{cases}
                \vvec(\coreop_{[i,j]}(M_{[i,j]})) = \mathfrak{X}(a_{i,j}, b_{i,j}, c_{i,j}, d_{i,j})  \vvec(M_{[i,j]}) & i, j \leq m,\\
                \vvec(\coreop_{[m+1,j]}(M_{[m+1,j]})) = \mathfrak{Y}(r_{j}, s_{j}) \vvec(M_{[m+1, j]}) & j \leq m,\\
                \coreop_{[m+1,m+1]}(M_{[m+1,m+1]}) = 0, &
            \end{cases}\\
            \text{where }&\lblines{
                a_{i,j} &= \sinc(\theta_i-\theta_j), &b_{i,j} = \cosc(\theta_i-\theta_j),\\
                c_{i,j} &= \sinc(\theta_i+\theta_j), &d_{i,j} = \cosc(\theta_i+\theta_j),\\
            }\text{ and } \lblines{
                r_j &= \sinc(\theta_j),\\
                s_j &= \cosc(\theta_j).\\
            }
        \end{aligned}
    \end{equation}
\end{proposition}
\begin{proof}
    As the formulae for the $2\times 2$ blocks and the $1\times 2$ blocks are essentially the same, the detailed derivation of the $2\times 2$ case is presented.

    Let $\xi := \vvec(M_{[i,j]}), \eta:= \vvec(\coreop_{[i,j]}(M_{[i,j]}))\in\realset^4$. Then, $M'_{[i,j]}=U^{\Her} M_{[i,j]}U$ is
    \[
        \begin{aligned}
            M'_{[i,j]} &= \frac{1}{2} \begin{bmatrix}
                \xi_1+\xi_4 - (\xi_2-\xi_3)\im & \xi_1-\xi_4 - (\xi_2+\xi_3)\im\\
                \xi_1-\xi_4 + (\xi_2+\xi_3)\im & \xi_1+\xi_4 + (\xi_2-\xi_3)\im\\
            \end{bmatrix}:=\frac{1}{2}\begin{bmatrix}
                \zeta_1 - \zeta_4\im & \zeta_2 - \zeta_3\im\\
                \zeta_2 + \zeta_3\im & \zeta_1 + \zeta_4\im
            \end{bmatrix}\\
        \end{aligned}
    \]
    The parameter block $\Phi_{[i,j]}, i, j \leq m$ is given by \eqref{eq:a2d-formulae}, which yields
    \begin{equation*}
        M''_{[i,j]} = \frac{1}{2}\begin{bmatrix}
            \zeta_1 -\im \zeta_4 & \zeta_2 -\im \zeta_3\\
            \zeta_2 +\im \zeta_3 & \zeta_1 +\im \zeta_4
        \end{bmatrix} \odot \begin{bmatrix}
            a - b\im & c - d\im\\
            c + d\im & a + b\im\\
        \end{bmatrix}.
    \end{equation*}
    Finally, $N_{[i,j]} = UM''_{i,j}U^{\Her}$ takes the form
    \begin{equation*}
        \begin{aligned}
            \begin{bmatrix}
                \eta_1 & \eta_3\\
                \eta_2 & \eta_4
            \end{bmatrix} &= U\left(\frac{1}{2}\begin{bmatrix}
                \zeta_1 -\im \zeta_4 & \zeta_2 -\im \zeta_3\\
                \zeta_2 +\im \zeta_3 & \zeta_1 +\im \zeta_4
            \end{bmatrix} \odot \begin{bmatrix}
                a - \im b & c - \im d\\
                c + \im d & a + \im b\\
            \end{bmatrix}\right)U^{\Her}\\
            &= \frac{1}{2}\begin{bmatrix}
                \zeta_1a + \zeta_2c - \zeta_3d - \zeta_4b & -\zeta_1b+ \zeta_2d + \zeta_3c - \zeta_4a\\
                \zeta_1b + \zeta_2d + \zeta_3c + \zeta_4a & \zeta_1a - \zeta_2c + \zeta_3d - \zeta_4b
            \end{bmatrix},\\
        \end{aligned}
    \end{equation*}
    which yields $\vvec(N_{[i,j]}) = \mathfrak{X}(a_{i,j}, b_{i,j}, c_{i,j}, d_{i,j})  \vvec(M_{[i,j]})$ for $i,j \leq m$ in~\eqref{eq:dskew-skew}.

    The formulae for $\coreop_{[m+1,j]}, j\leq m$ are obtained in the same but simpler procedure. For the left-over diagonal $\coreop_{[m+1,m+1]}$, it suffices to notice that an $1\times 1$ skew-symmetric matrix is always $0$.
\end{proof}

Note that the formulae for the $2\times 1$ blocks $N_{[i,m+1]}, i \leq m$ are not necessary due to the skew-symmetry in $N$. In fact, all upper triangular blocks $N_{[i,j]}, i < j$ can be obtained from the lower counterpart via $N_{[i,j]}= -N_{[j,i]}^{\T}$. Furthermore, \cref{coro:core-diag-blocks} shows that the diagonal blocks $N_{[i,i]}$ are equal to $M_{[i,i]}$.

\begin{corollary}\label{coro:core-diag-blocks}
    For all $M_{[i,i]}\in\skewm_2$ and all $\theta_i\in\realset$ that determines $\coreop_{[i,i]}$, it holds that $\coreop_{[i,i]}(M_{[i,i]}) = M_{[i,i]}$.
\end{corollary}
\begin{proof}
    The core action at the $i$th diagonal block is given by
    \begin{equation*}
        \mathfrak{X}(1, 0, c, d) = \frac{1}{2}\sbmatrix{
            1+c & -d & -d & 1-c\\
            d & 1-c &-1-c & -d\\
            d & -1+c & 1+c &-d\\
            1-c & d & d & 1+c\\
        }
    \end{equation*}
    as $a = \sinc(\theta_i-\theta_i) = 1$ and $b = \cosc(\theta_i-\theta_i)=0$. The values of $c=\sinc(2\theta_i)$ and $d=\cosc(2\theta_i)$ do not matter as they are canceled by $\skewm_2\ni M_{[i,i]} = \sbmatrix{
        0 & -x\\
        x & 0
    }$:
    \begin{equation*}
        \frac{1}{2}\sbmatrix{
            1+c & -d & -d & 1-c\\
            d & 1-c &-1-c & -d\\
            d & -1+c & 1+c &-d\\
            1-c & d & d & 1+c\\
        }\sbmatrix{
            0\\
            x\\
            -x\\
            0
        } = \sbmatrix{
            0\\
            x\\
            -x\\
            0
        },\forall c,d\in\realset,
    \end{equation*}
    i.e., $\coreop_{[i,i]}(M_{[i,i]}) = M_{[i,i]},\forall M_{[i,i]}\in\skewm_2$ holds for all $\coreop_{[i,i]}$ determined by any $\theta_i\in\realset$.
\end{proof}

\subsection{Inverse of the Core Map}

Recall that, for all $X\in\skewm_n$, the differentiation of the exponential is given by $\dexpof{A}[X] =Q  \dskew_A(X) = Q R  \core(R^{\T}XR) R^{\T}$.
If the differentiation is invertible, then its inverse can be expressed in terms of $\dskew_A^{-1}$:
\begin{equation}\label{eq:dexp-skew-inv}
    \begin{cases}
        \dexpinvof{A}[\Delta] = \dskew_A^{-1}\left(Q^{\T} \Delta\right) &\forall \Delta \in T_Q\so_n,\\
        \dskew_A^{-1}(Y) = R  \core^{-1}(R^{\T}YR) R^{\T} & \forall Y\in\skewm_n .\\
    \end{cases}
\end{equation}
These expressions imply that the invertibility of $\dexpof{A}$ and the linear transformations on $\skewm_n$ are equivalent as given in \cref{prop:invertible-equivalence}.

\begin{proposition}\label{prop:invertible-equivalence}
    For $A\in\skewm_n$, if any one of the linear maps: (i) $\dexpof{A}:\skewm_n\to T_Q\so_n$; (ii) $\dskew_A:\skewm_n\to\skewm_n$; and (iii) $\core:\skewm_n\to \skewm_n$ is invertible, the other two linear maps are also invertible.
\end{proposition}
\begin{proof}
    The equivalence between the invertibility follows from the fact that the left multiplication by $Q$ is a bijection between $\skewm_n$ and $T_Q\so_n$, and the congruence by orthogonal Schur vectors $R$ is a bijective linear transformation of $\skewm_n$.
\end{proof}

Since $\dexpinvof{A}:T_Q\so_n\to \skewm_n$ is well-defined if and only if $\core^{-1}$ is well-defined, the goal is to determine the invertibility of $\core$ and find the symbolic inverse formula. The first step is to obtain the symbolic inversions of the matrices in~\eqref{eq:dexp-skew-mat}.

\begin{lemma}\label{lemma:dexp-skew-inv-mat}
    The matrices $\mathfrak{X}(a,b,c,d)$ and $\mathfrak{Y}(r,s)$ in~\eqref{eq:dexp-skew-mat} are invertible if and only if the pairs $(a,b)$, $(c,d)$ and $(r,s)$ are not $(0,0)$, and their inverses are given by
    \begin{equation}\label{eq:dexp-skew-inv-mat}
        \lblines{
            &\mathfrak{X}^{-1}(a,b,c,d) = \mathfrak{X}\left(\frac{a}{a^2+b^2}, \frac{-b}{a^2+b^2}, \frac{c}{c^2+d^2}, \frac{-d}{c^2+d^2}\right):= \mathfrak{X}\left(a',b',c',d'\right),\\
            &\mathfrak{Y}^{-1}(r,s) = \mathfrak{Y}\left(\frac{r}{r^2+s^2}, \frac{-s}{r^2+s^2} \right):=\mathfrak{Y}(r',s').\\
        }
    \end{equation}
\end{lemma}
\begin{proof}
    It is readily checked that $\mathfrak{X}\left(a',b',c',d'\right) \mathfrak{X}(a,b,c,d) = I_4$ when $a', \ldots, d'$ are well-defined, and that $\mathfrak{Y}\left(r',s'\right)  \mathfrak{Y}\left(r,s\right) = I_2$ when $r', s'$ are well-defined. Otherwise, the corresponding matrix is not invertible.
\end{proof}

Recall that the actions of $\core$ at the diagonal blocks are the identity, which makes the inverse actions of $\core$ of the diagonal blocks the identity. Then, applying the symbolic inversions~\eqref{eq:dexp-skew-inv-mat} on the off-diagonal blocks yields the symbolic inversion.

\begin{proposition}\label{prop:core-inverse}
    If $\core:\skewm_n\to\skewm_n$ is invertible, its inversion acting on all $N\in\skewm_n$ by blocks $N_{[i,j]}$ with linear maps $\coreop_{[i,j]}^{-1}, i,j\leq k$ is given by
    \begin{equation}\label{eq:core-inverse}
        \sbmatrix{
                M_{[1,1]} & \cdots & M_{[1,k]}\\
                \vdots & \ddots & \vdots\\
                M_{[k,1]} & \cdots & M_{[k,k]}
        }
        = M = \core^{-1}(N) =
        \sbmatrix{
            \coreop_{[1,1]}^{-1}(N_{[1,1]}) & \cdots & \coreop_{[1, k]}^{-1}(N_{[1,k]})\\
            \vdots & \ddots & \vdots\\
            \coreop_{[k,1]}^{-1}(N_{[k,1]}) & \cdots & \coreop_{[k, k]}^{-1}(N_{[k,k]})\\
        }.
    \end{equation}
    For the vectorized $\vvec(N_{[i,j]})\in\realset^4, i,j\leq m$ and $\vvec(N_{[m+1,j]})\in\realset^2, j\leq m$ for $n = 2m+1$, the linear maps are given by
    \begin{equation}\label{eq:core-vec-inverse}
        \begin{aligned}
            &\begin{cases}
                \vvec(\coreop_{[i,i]}^{-1}(N_{[i,i]})) = \vvec(N_{[i,i]}) & i\leq k,\\
                \vvec(\coreop_{[i,j]}^{-1}(N_{[i,j]})) = \mathfrak{X}(a'_{i,j}, b'_{i,j}, c'_{i,j}, d'_{i,j})  \vvec(N_{[i,j]}) & j< i \leq m,\\
                \vvec(\coreop_{[m+1,j]}^{-1}(N_{[m+1,j]})) = \mathfrak{Y}(r'_{j}, s'_{j}) \vvec(N_{[m+1, j]}) & j \leq m,\\
                \coreop_{[m+1,m+1]}^{-1}(N_{[m+1,m+1]}) = 0, &
            \end{cases}\text{ where }\\
            &\lblines{
                a'_{i,j} &= \cotc\left((\theta_i - \theta_j)\big\slash 2\right), &b'_{i,j} = -(\theta_i - \theta_j)\big\slash 2,\\
                c'_{i,j} &= \cotc\left((\theta_i + \theta_j)\big\slash 2\right), &d'_{i,j} = -(\theta_i + \theta_j)\big\slash 2,\\
            }\text { and } \lblines{
                r'_j &= \cotc\left(\theta_j\big\slash 2\right),\\
                s'_j &= -\theta_j\big\slash 2.\\
            }
        \end{aligned}
    \end{equation}
\end{proposition}
\begin{proof}
    Substituting~\eqref{eq:dexp-skew-inv-mat} into the inversion of~\eqref{eq:dskew-skew} yields the vectorized forms of the independent linear maps in~\eqref{eq:core-vec-inverse}. Notice the trigonometric identity of the squared sum
    $\sinc(z)^2 + \cosc(z)^2 = (2 - 2\cos(z)) \big/ z^2 = 4\sin^2(z/2)\big/ z^2$. Substitute it to the fractions $\sinc(z)\big/(\sinc(z)^2 + \cosc(z)^2)$ and $-\cosc(z)\big/(\sinc(z)^2 + \cosc(z)^2)$ to get
    \begin{equation*}
        \lblines{
            &\frac{\sin(z)}{z} \bigg/ \frac{4\sin(z/2)^2}{z^2} = \frac{2\sin(z/2)\cos(z/2)}{z} \bigg/ \frac{4\sin(z/2)^2}{z^2} = \frac{\cot(z/2)}{z/2} = \cotc\left(\frac{z}{2}\right),\\
            &-\frac{\cos(z)-1}{z} \bigg/ \frac{2-2\cos(z)}{z^2} = -\frac{z}{2}.
        }
    \end{equation*}
    Therefore, $a', \ldots, s'$ are given by the substitutions $z=\theta_i-\theta_j$, $\theta_i+\theta_j$ and $\theta_j$.
\end{proof}

It is important to include the extra specification on the diagonal blocks in~\eqref{eq:core-vec-inverse}, compared to~\eqref{eq:dskew-skew}, as it is possible to have $\pm2\theta_i = 2l\pi$ for some $i\leq m$ and $l\neq 0$ on the diagonal block, where the associated matrix $\mathfrak{X}(1,0,0,d)$ is rank-deficient but the inverse action remains the well-defined identity. The identity action of the diagonal blocks yields a necessary and sufficient condition of invertibility in \cref{thm:invertibility-skew-symm}, which is different from the invertibility of $\dexpof{A}:\realmat{n}\to T_Q\gl_n$ for $Q = \exp(A)$.

\begin{theorem}\label{thm:invertibility-skew-symm}
    For $A\in\skewm_n$ with the angles $\Theta\in\realset^m$, $\core:\skewm_n\to\skewm_n$ and $\dexpof{A}:\skewm_n\to T_Q\so_n$ are invertible if and only if
    \begin{equation}\label{eq:invertible-condition}
        \lblines{
            \theta_i\pm\theta_j\neq 2l\pi, &\forall j < i \leq m , l=\pm1,\pm2,\ldots\\
            \theta_j\neq 2l\pi, &\forall j\leq m, l=\pm1,\pm2,\ldots\text{ when } n=2m+1.
        }
    \end{equation}
\end{theorem}
\begin{proof}
    Since the entries in the off-diagonal blocks $N_{[i,j]}, i\neq j$ of a skew-symmetric matrix $N$ are free, the linear map $\coreop^{-1}_{[i,j]}$ is invertible if and only if the respective $\mathfrak{X}$- or $\mathfrak{Y}$-matrix is invertible. Therefore, $\core$ is invertible if and only if all $\mathfrak{X}$- and $\mathfrak{Y}$-matrices involved with the off-diagonal $N_{[i,j]}, i\neq j$ have their symbolic inversions well-defined. According to \cref{eq:dexp-skew-mat} and \cref{eq:dexp-skew-inv-mat}, this occurs if and only if $z\neq 2l\pi, \forall l=\pm1,\pm2,\ldots$, where $z$ takes values from $\theta_i\pm\theta_j, j < i\leq m$ and $\theta_j, j\leq m$.
\end{proof}

Note that the classic invertibility condition~\cite[Prop.~7, Sec.~1.2]{rossmann2006lie} states that $\dexpof{A}:\realmat{n} \to T_Q\gl_n$ with $A\in\realmat{n}$ is invertible if and only if for all pair of eigenvalues $\lambda_i, \lambda_j, \forall i,j\leq n$, it holds that $\lambda_i-\lambda_j \neq 2l\pi\im, l = \pm1, \ldots$. For $A\in\skewm_n$ with the angles $\Theta$, the eigenvalues are either $\pm\theta_i\im, 1\leq i \leq m$ or $0$. Then,  $\im\theta_i\pm \im\theta_j = 2l\pi\im$ is equivalent with $\theta_i\pm \theta_j = 2l\pi$. That is to say $\dexpof{A}:\realmat{n}\to T_Q\gl_n, A\in\skewm_n$ is invertible if and only if
\begin{equation}\label{eq:inv-necessary}
    \lblines{
        \theta_j\pm\theta_i\neq 2l\pi, &\forall i , j \leq m , \forall l= \pm1, \pm2,\ldots \\
        \theta_j\neq 2l\pi, &\forall j\leq m, \forall l= \pm1, \pm2,\ldots \text{ when } n=2m+1.
    }
\end{equation}
Observe that condition~\eqref{eq:inv-necessary} is stronger than condition~\eqref{eq:invertible-condition}, where the condition on the $i=j$ case in relaxed. This is due to the fact that the domain and co-domain considered for $\dexpof{A}$ are larger for~\eqref{eq:inv-necessary} than for~\eqref{eq:invertible-condition}.

\subsection{Algorithm and Complexity}

The pseudocodes for computing $\dskew_A(X)$ in $\dexpof{A}[X] = Q \dskew_A(X)$, \cref{alg:differentiation}, and its inversion $\dskew^{-1}_A(Y)$, \cref{alg:inversion-differentiation}, are presented along with computational cost respectively in this section.

\subsubsection{Differentiation of the Matrix Exponential}

For the differentiation $\dexpof{A}[X] = Q \dskew_A(X), \forall X\in\skewm_n$, \cref{alg:differentiation} generates the skew-symmetric output $\dskew_A(X) = R\left(\core(R^{\T}XR)\right)R^{\T}$ given the Schur decomposition $(\Theta, \Xi, R)$ of $A$. When the full differentiation is needed, additional work is needed to compute $A\mapsto R\exp(\Xi )R^{\T} = Q$ and $\dskew_A(X) =: Y \mapsto Q Y$. It is convenient to have direct access to $Y\in\skewm_n$, which characterizes a tangent vector $Q Y\in T_Q\so_n$ intrinsically. For example, the geodesic in $\so_n$ that emanates from $Q$ along $Q Y$ is given by $Q \exp(t Y), t\in [0,\infty)$, c.f.~\cite{edelman1998}. Then, working with $\dskew_A(X)$ directly is more convenient for geodesics-related applications on $\so_n$. For these reasons, \cref{alg:differentiation} presents the computation of $\dskew_A$ as a self-contained routine rather than the full differentiation $\dexpof{A}[X] = Q \dskew_A(X)$. The routine consists of two stages, the preprocessing stage from \cref{dexp-alg:preprocessing-start} to \cref{dexp-alg:preprocessing-end} and the computation stage from \cref{dexp-alg:computing-start} to \cref{dexp-alg:computing-end}. When there are multiple calls to \cref{alg:differentiation} at the same $A$ with different perturbations $X$, the preprocessing stage is executed once and matrices/parameters $\mathfrak{X}_{i,j}, \mathfrak{Y}_j$ are stored for the multiple calls to $\dskew_A(X)$.

\begin{algorithm}

    \caption{Linear map $\dskew_A:\skewm_n\to\skewm_n$ in $\dexpof{A}[X]=Q \dskew_A(X)$}\label{alg:differentiation}
    \begin{algorithmic}[1]
        \REQUIRE Schur decomposition $(\Theta, \Xi, R)$ of $A$.
        \FOR{$j = 1,\ldots, k$}\label{dexp-alg:preprocessing-start}
            \FOR{$i = j+1,\ldots, m$}
                \STATE $a_{i,j}\gets \sinc(\theta_i-\theta_j)$ and $b_{i,j}\gets \cosc(\theta_i-\theta_j)$
                \STATE $c_{i,j}\gets \sinc(\theta_i+\theta_j)$ and $d_{i,j}\gets \cosc(\theta_i+\theta_j)$
                \STATE $\mathfrak{X}_{i,j}\gets \mathfrak{X}(a_{i,j}, b_{i,j}, c_{i,j}, d_{i,j})$ given by~\eqref{eq:dexp-skew-mat} 
            \ENDFOR
            \IF{$n==2m+1$}
                \STATE $r_{j}\gets \sinc(\theta_j)$ and $s_{j}\gets \cosc(\theta_j)$
                \STATE $\mathfrak{Y}_{j}\gets \mathfrak{Y}(r_j,s_j)$ given by~\eqref{eq:dexp-skew-mat} 
            \ENDIF
        \ENDFOR\label{dexp-alg:preprocessing-end}\codecomment{Preprocessing done.}
        \STATE\label{dexp-alg:computing-start} $M\gets R^{\T}X R$ \codecomment{Start Computing.}
        \FOR{$j = 1,\ldots, k$}
            \STATE $N_{[j, j]}\gets M_{[j,j]}$
            \FOR{$i = j+1,\ldots, m$}
                \STATE $\vvec(N_{[i,j]})\gets \mathfrak{X}_{i,j} \vvec(M_{[i,j]})$
                \STATE $N_{[j,i]} \gets -N_{[i,j]}^{\T}$
            \ENDFOR
            \IF{$n==2m+1$}
                \STATE $\vvec(N_{[m+1,j]})\gets \mathfrak{Y}_{j} \vvec(M_{[m+1,j]})$\
                \STATE $N_{[j,m+1]} \gets -N_{[m+1,j]}^{\T}$
            \ENDIF
        \ENDFOR
        \RETURN \label{dexp-alg:computing-end} $Y\gets R N R^{\T}$
    \end{algorithmic}
\end{algorithm}

\cref{tab:dexp-ops} reports the complexity analysis on the number of operations (ops in short) on floating point data. In particular, the operation required in the proposed formula consists of (i) The Schur decomposition of $A$ proposed in~\cite{mataigne2024eigenvalue}, which requires $14/3n^3$ ops; (ii) the execution of \cref{alg:differentiation}, which requires $6n^3$ ops for the four matrix multiplications in $M\gets R^{\T}XR$ and $Y\gets RNR^{\T}$; (iii) (Optional) For a full differentiation $\dexpof{A}[X] = Q \dskew_A(X)$, an extra exponential and matrix multiplication are performed, which require $4n^3$ ops. Note that the skew-symmetry of $X$ and $Y$ saves $2n^3$ ops at step (ii), and that the Schur decomposition of $A$ is reused at step (iii) such that the exponential only requires $2n^3$ ops.

Also note that the total ops in the Schur decomposition and matrix exponential is $20/3n^3$ while the scaling and squaring algorithm~\cite{higham2005scaling} requires $20/3 n^3$ ops for the $[3/3]$ Pad{\'e} approximant, or $(46/3+4\sigma)n^3$ ops the $[13/13]$ Pad{\'e} approximant, i.e., the exponential computed by the Schur decomposition of skew-symmetric matrices is no slower than the scaling and squaring algorithm. The numerical experiments in \Cref{sec:experiments}, see \cref{sub@fig:expm-time}, support these observations on the complexity.

The complex formulae~\eqref{eq:dexp-diag} and~\eqref{eq:linear-diag} compute the skew-symmetric output $X\mapsto \dskew_A(X)$ and the full differential output $X\mapsto \dexpof{A}[X]$, respectively. Note that the eigendecomposition of non-symmetric matrices is significantly slower than the specialized Schur decomposition of skew-symmetric matrices proposed in~\cite{mataigne2024eigenvalue}. Therefore, the required eigendecomposition of $A\in\skewm_n$ is obtained from the Schur decomposition via~\eqref{eq:schur-to-spectral} with an extra $O(n^2)$ ops. On the other hand, the Pad{\'e} approximant proposed in~\cite{al2009exponential} demands the Pad{\'e} approximant $\poly{r/r}{\exp}(A)$ with $r = 3,5,7,9,13$. The Pad{\'e} approximants consists of a few matrix multiplications, one LU-decomposition and two calls to the LU-solver. Depending on the Pad{\'e} order $r$ and the scaling order $\sigma$, the number of matrix multiplications varies.

\begin{table}[!tbp]
    \centering
    \caption{Complexity of $\dexpof{A}[X]$ and $\dskew_{A}(X)$, $\forall A,X\in\skewm_n$}
    \begin{tabular}{|l|c|c:c|}
        \hline
        Formula & Input: Exponential or & \multicolumn{2}{c|}{Output Format}\\
         & Decomposition of $A$  & $\dskew_A(X)$& $\dexpof{A}[X]$ \\
        \hline
        Pad{\'e} $[3/3]$$^{\dagger}$ & $(6 + 2\big/ 3)n^3$ for $\exp(A)$ & $12n^3$ & $10n^3$ \\
        \hline
        Pad{\'e} $[13/13]$$^{\dagger}$ & $(14 + 2\sigma+ 2\big/ 3)n^3$ for $\exp(A)$ & $(30 + 4\sigma)n^3$ & $(28+4\sigma)n^3$ \\
        \hline
        Complex$^*$ & $14\big/3 n^3$ for $A=V\Lambda V^{-1}$ & $24n^3$ in~\eqref{eq:linear-diag} & $28n^3$ in~\eqref{eq:dexp-diag}\\
        \hline
        Real$^\#$ & $14\big/3 n^3$ for $A=R\Xi R^{\T}$ & $6n^3$ & $10 n^3$\\
        \hline
		\multicolumn{4}{l}{\tiny $\dagger$~\cite{al2009exponential}: The Pad{\'e} approximant \:\:  $*$~\cite{najfeld1995}: Daletski$\breve{\mathrm{\i}}$--Kre$\breve{\mathrm{\i}}$n Formulae~\eqref{eq:linear-diag} and~\eqref{eq:dexp-diag} \:\: $\#$: The proposed formula~\eqref{eq:dskew-skew}.}\\
    \end{tabular}
    \label{tab:dexp-ops}
\end{table}

\subsubsection{Inverse of the Differentiation}

For the inverse of the differentiation $\dexpinvof{A}[Q Y] = \dskew_A^{-1}(Y), \forall Y\in\skewm_n$, \cref{alg:inversion-differentiation} takes a skew-symmetric input $Y\in\skewm_n$ and generates $\dexpinvof{A}[Q Y]$. When the full differentiation $\Delta=QY \in T_Q\so_n$ is provided, the extra computations for $A\mapsto \exp(A) = Q$ and $\Delta \mapsto Q^{\T} \Delta$ are required. \Cref{alg:inversion-differentiation} has almost the same structure with \cref{alg:differentiation} except for the extra checks of the invertibility. The discussion of complexity presented in \cref{tab:dlog-ops} is also similar to \cref{tab:dexp-ops} except for the same $O(n^3)$ ops in the first and second rows of different Pad{\'e} approximant. According to~\cite{al2013logarithm}, different orders in the Pad{\'e} approximant for the differentiation of the matrix logarithm only affects the number of system solves in the block diagonal $\Xi$ in $A=R\Xi R^{\T}$, which only requires $O(n^2)$ ops. Therefore, the order does not affect the $O(n^3)$ term for $A\in\skewm_n$. Note that the nearly identical elapsed times of the $[1/1]$ and $[7/7]$ Pad{\'e} approximants reported in \Cref{fig:timing_dlog_skew} and \Cref{fig:timing_dlog_full} support the claim that they have the same computational complexity.

\begin{algorithm}
    \caption{Linear map $\dskew_A^{-1}$ in $\dexpinvof{A}[Q Y]= R  \core^{-1}(R^{\T}YR) R^{\T}$.}\label{alg:inversion-differentiation}
    \begin{algorithmic}[1]
        \REQUIRE Schur decomposition $(\Theta, \Xi, R)$ of $A$.
        \FOR{$j = 1,\ldots, k$}\label{dlog-alg:preprocessing-start}
            \FOR{$i = j+1,\ldots, m$}
                \STATE $a'_{i,j}\gets \cotc((\theta_i-\theta_j)\big/ 2)$ and $b'_{i,j}\gets -(\theta_i-\theta_j)\big/ 2$
                \STATE $c'_{i,j}\gets \cotc((\theta_i+\theta_j)\big/ 2)$ and $d'_{i,j}\gets -(\theta_i+\theta_j)\big/ 2$
                \STATE Report rank-deficient warnings if $\exists l\neq 0, s.t., b'_{i,j} \text{ or } d'_{i,j} = l\pi$\codecomment{Sanity check}
                \STATE $\mathfrak{X}_{i,j}\gets \mathfrak{X}(a'_{i,j}, b'_{i,j}, c'_{i,j}, d'_{i,j})$ given by~\eqref{eq:dexp-skew-mat}
            \ENDFOR
            \IF{$n==2m+1$}
                \STATE $r'_{j}\gets \sinc(\theta_j\big/ 2)$ and $s'_{j}\gets -\theta_j\big/ 2$
                \STATE Report rank-deficient warnings if $\exists l\neq 0, s.t., s'_j = l\pi$ \codecomment{Sanity check}
                \STATE $\mathfrak{Y}_{j}\gets \mathfrak{Y}(r'_j,s'_j)$ given by~\eqref{eq:dexp-skew-mat}
            \ENDIF
        \ENDFOR\label{dlog-alg:preprocessing-end}\codecomment{Preprocessing done.}
        \STATE Execute \cref{dexp-alg:computing-start} to \cref{dexp-alg:computing-end} in \cref{alg:differentiation}. \codecomment{Computing Stage.}
    \end{algorithmic}
\end{algorithm}

\begin{table}[!tbp]
    \centering
    \caption{Complexity of $\dexpinvof{A}[Q Y]$ and $\dskew_A^{-1}(Y)$, $\forall A,Y\in\skewm_n$}
    \label{tab:dlog-ops}
    \begin{tabular}{|l|c|c:c|}
        \hline
        Formula & Input Parameters:  & \multicolumn{2}{c|}{Input Format}\\
         & Decomposition of $A$  & $Y\in\skewm_n$& $Q Y\in T_Q\so_n$ \\
        \hline
        Pad{\'e} $[1/1]$$^{\dagger}$ & $(4 + 13\big/ 6)n^3$ for $A = R\Xi R^{\T}$ & $(24+2\sigma)n^3$ & $(22+2\sigma)n^3$ \\
        \hline
        Pad{\'e} $[7/7]$$^{\dagger}$ & $(4 + 13\big/ 6)n^3$ for $A = R\Xi R^{\T}$ & $(24+2\sigma)n^3$ & $(22+2\sigma)n^3$ \\
        \hline
        Complex$^*$ & $14\big/3 n^3$ for $A=V\Lambda V^{-1}$ & $24n^3$ in~\eqref{eq:linear-diag} & $28n^3$ in~\eqref{eq:dexp-diag}\\
        \hline
        Real$^\#$ & $14\big/3 n^3$ for $A=R\Xi R^{\T}$ & $6n^3$ & $10 n^3$\\
        \hline
		\multicolumn{4}{l}{\tiny $\dagger$~\cite{al2013logarithm}: The Pad{\'e} approximant \:\:  $*$~\cite{najfeld1995}: Daletski$\breve{\mathrm{\i}}$--Kre$\breve{\mathrm{\i}}$n Formulae~\eqref{eq:linear-diag} and~\eqref{eq:dexp-diag} \:\: $\#$: The proposed formula~\eqref{eq:dexp-skew-inv}.}\\
    \end{tabular}

\end{table}

\section{Nearby Matrix Logarithm: A Well-Defined Local Inversion}\label{sec:nearby-logarithm}

This section further investigates the set of skew-symmetric matrices with rank-deficient differentiation on the matrix exponential. The nearby matrix logarithm is developed as a well-defined local inversion to the matrix exponential on $\so_n$ such that it works beyond the image of the principal logarithm in $\skewm_n$. Such a tool is important for geometric analysis on manifolds with special orthogonally constraints, e.g.,~\cite{jurdjevic2020extremal}.  Preliminary experiments conducted in this section demonstrate the ability of the nearby matrix logarithm to be continuous onto a subset of $\skewm_n$ that is not included in the image of the principal logarithm.

\subsection{Further Characterization on the Invertibility}

The invertibility of $\dexpof{A}$ is further investigated in this section by the angular characterization $\Theta\in\realset$ of $A\in\skewm_n$ and the matrix $2$-norm $\skewm_n$, where $\|A\|_2 = \max_{i\leq m}|\theta_i|$. \Cref{lemma:angular-deviation} reveals the relationship between the angles and the matrix $2$-norm via Weyl's inequality~\cite{horn2012matrix} on the eigenvalues of the Hermitian matrices $\im A$ and $\im B$. Recall that these eigenvalues are associated with the angles of $A$ and $B$ via~\eqref{eq:schur-to-spectral}.

\begin{lemma}\label{lemma:angular-deviation}
    For $A\in\skewm_n$ with angles $\Theta\in\realset^m$ and $C\in\skewm_n$ with $\|C\|_2 < \delta$, there exists a vector of angles $\Sigma=(\sigma_1,\ldots, \sigma_m)\in\realset^m$ of $B := A+C$ that satisfies
    \begin{equation}\label{eq:angular-deviation}
        |\sigma_i - \theta_i| < \delta, \forall i = 1,\ldots, m.
    \end{equation}
\end{lemma}
\begin{proof}
    Suppose $\Gamma\in\realset^m$ is a set of angles of $C$. Since $\|C\|_2 = \max_i|\gamma_i|$, it holds that $|\gamma_i| < \delta, \forall i = 1,\ldots,m$.
    Notice that $\im A$ and $\im C$ are Hermitian with the eigenvalues $\{\theta_1, -\theta_1,\ldots\}$ and $\{\gamma_1,-\gamma_1, \ldots\}$. Then, Weyl's inequality bounds the deviations of eigenvalues in $\im B$ from $\im A$ as $|\lambda_i(\im B) - \lambda_i(\im A)| \leq \max_{j}|\lambda_j(\im C)| < \delta$
    where $\lambda_i$ returns the respective eigenvalue in an appropriate order. Since $\lambda_i(\im B)$ are the signed angles in $\Sigma$, the bounds~\eqref{eq:angular-deviation} on the corresponding angles follow.
\end{proof}

\begin{definition}\label{def:conj-locus}
    The set of $A\in\skewm_n$ where $\dexpof{A}:\skewm_n\to\so_n$ is rank-deficient is termed the \emph{tangent conjugate locus}.\footnote{
    The term \emph{tangent conjugate locus} is adequate in the following sense. On the Lie group $\so_n$, consider the Frobenius Riemannian metric (but note that, in this paper, ``$\mathrm{dist}$'' stands for the distance in the spectral norm). Since the metric is bi-invariant~\cite[Example~4.10]{GuiguiMiolanePennec2023}, the matrix exponential (which is the group exponential) coincides with the Riemannian exponential $\mathrm{Exp}_{I_n}$; see, e.g.,~\cite[Proposition~21.20(3)]{GallierQuaintance2020}. By definition, the \emph{tangent conjugate locus} in $T_{I_n}\so_n$ is the set of initial velocities $v\in T_{I_n}\so_n$ where $\der \Exp_{I_n}(v): T_{v}T_{I_n}\so_n\to T_{\Exp_{I_n}(v)}\so_n$ is rank-deficient.  Note that $T_{I_n}\so_n = \skewm_n$, $T_A\skewm_n = \skewm_n$ and $\Exp_{I_n}(A) = \exp(A)$ hold for all $A\in\skewm_n$.} In view of \cref{thm:invertibility-skew-symm}, it is
    \begin{equation}\label{eq:conjugate-locus}
        \fraks:= \left\{
                A\in\skewm_n:
                \begin{aligned}&\exists i\neq j, \mathbb{Z}\ni l\neq 0, \text{ s.t. } \theta_i\pm \theta_j = 2l\pi, \\
                &\text{ or s.t. }\theta_i = 2l\pi \text{ when } n = 2m+1
            \end{aligned}
        \right\}.
    \end{equation}
\end{definition}

The tangent conjugate locus $\fraks$ is described by a countable collection of conditions indexed by the integers $l > 0$, where the $l<0$ cases are subsumed by flipping the angles. Unless otherwise specified, $l$ in the rest of this paper is always a positive integer. The countably many subsets are denoted as
\begin{equation}\label{eq:closed-subset}
    \lblines{
        &\conjset_{l,\pm} :=  \{A\in\skewm_n:\exists \theta_i \pm \theta_j=2l\pi\},\\
        &\conjset_{l,*} := \{A\in\skewm_n:\exists\theta_i=2l\pi\},
    }
\end{equation}
such that $\fraks = \bigcup_{l > 0}\conjset_{l,\pm}$ when $n=2m$, and $\fraks = \bigcup_{l >0}\left(\conjset_{l,\pm} \cup \conjset_{l,*}\right)$ when $n = 2m+1$.

The matrix exponential $\exp:\skewm_n\to\so_n, n\geq 3$ is not an open map as shown in~\cite[Proposition 8.3]{absil2025ultimate}. \Cref{prop:loss-of-open-map} presents a necessary condition for the non-openness in $\exp:\skewm_n\to \so_n$, whose contrapositive statement can be used to extract certain subsets of $\skewm_n$ where $\exp$ is an open map.

\begin{proposition}\label{prop:loss-of-open-map}
    If the matrix exponential $\exp:\mathcal{X}\to \mathcal{Y}$, with $\mathcal{X}$ an open subset of $\skewm_n$, is not an open map, then $\mathcal{X}\cap \fraks \neq \emptyset$, and consequently there exist $Q\in\mathcal{Y}$ with repeated eigenvalues.
\end{proposition}
\begin{proof}
    Since the matrix exponential $\exp:\skewm_n\to\so_n$ is a smooth map, it is necessary (but not sufficient) to have rank-deficient $\dexpof{A}:\skewm_n\to T_Q\so_n$ at some $A\in\mathcal{X}$ if $\exp:\mathcal{X}\to \mathcal{Y}$ is not an open map, i.e., $\mathcal{X}\cap \fraks \neq \emptyset$. For such an $A\in\fraks$ with Schur decomposition $(\Theta, \Xi, R)$, at least one of the following conditions holds for some integer $0\neq l\in\mathbb{Z}$: (i) $\exists i\neq j,\theta_i = \theta_j + 2l\pi$ (ii) $\exists \theta_i = 2l\pi$ and $n = 2m+1$. Consider $Q = \exp(A) = R\exp(\Xi ) R^{\T}$. In (i), the block $\exp(\Xi_{[i,i]}) = \exp(\Xi_{[j,j]}), i\neq j$ yields repeated complex eigenvalues of $\exp(\pm\theta_i\im)=\exp(\pm\theta_j\im)$ with at least multiplicities of two. In (ii), $\exp(\Xi_{[i,i]}) = I_2$ and the left-over block $\exp(\Xi_{[k,k]}) = 1$ yields a repeated real eigenvalue of one with at least a multiplicity of three.
\end{proof}

\Cref{coro:differential-skew-open-subsets} shows that  $\skewm_n \setminus \fraks$ (i.e., the set of $A$ where $\dexpof{A}$ is invertible) is an open set with countably many open and connected components.

\begin{corollary}\label{coro:differential-skew-open-subsets}
    The tangent conjugate locus $\fraks$ is a closed set in $\skewm_n$ and $\skewf:=\skewm_n\setminus \fraks$ is an open set\footnote{$\skewf$ is not an open set in $\realmat{n}$, since $\skewm_n$ is a subspace of $\realmat{n}$ but not an open set.} in $\skewm_n$. Furthermore, $\skewf$ consists of countably many connected open components, denoted by a countable set of indices $\mathcal{E}$ as follows:
    \begin{equation}\label{eq:differential-skew-open-subsets}
        \bigcup_{e\in \mathcal{E}} \skewf_{e} := \skewf=\skewm_n\setminus \fraks.
    \end{equation}
\end{corollary}
\begin{proof}
    Notice that $\fraks$ is a union of countably many closed subsets $\conjset_{l,\pm}$, (and $\conjset_{l,*}$ if $n = 2m+1$). Recall that the union of a locally finite collection of closed subsets remains closed, where a collection of subsets is locally finite if, for all points in the space, there exists a neighborhood that only intersects finitely many subsets. For $A,B\in\skewm_n$ where $\|A-B\|_2<\pi/2$, in view of \cref{lemma:angular-deviation}, there exist the respective angles $\Theta, \Sigma\in\realset^m$ such that $|(\theta_i\pm \theta_j) - (\sigma_i\pm \sigma_j)| < \pi$, for all $i,j\leq m$, i.e., $\sigma_i\pm \sigma_j$ takes its value from the interval $(\theta_i\pm\theta_j - \pi, \theta_i\pm\theta_j + \pi)$, for all $i,j \leq m$. For each interval with length $2\pi$, there is at most one $2l\pi\in (\theta_i\pm\theta_j - \pi, \theta_i\pm\theta_j + \pi), l\in\mathbb{Z}$. In total, there are at most $2m$ conditions $\sigma_i\pm\sigma_j = 2l\pi$ (with different $l$) that can be met in all balls $\{B:\in\skewm_n:\|A-B\|_2\leq \varepsilon < \pi/2\}$ centered at any $A\in\skewm_n$, i.e., the ball with radius $\varepsilon < \pi / 2$ intersects at most $2m$ different $\conjset_{l,\pm}$. Similarly, any ball with radius $\varepsilon < \pi$ intersects at most $m$ different $\conjset_{l, *}$. Therefore, the collection of the sets in~\eqref{eq:closed-subset} is locally finite, hence its union $\fraks$ is closed and $\skewm_n\setminus \fraks$ is open.

    Furthermore, $\skewf$ has at most countably many open and connected components as it is an open set in a (second-countable) Euclidean space $\skewm_n$. Finally, $\skewf$ has infinitely many connected components. Indeed, for $A$ with $\|A\|_2 < l\pi, l\in\mathbb{Z}$, there exists a component beyond the ball $\{A\in\skewm_n:\|A\|_2 > 3l\pi\}$ separated by $\conjset_{l,\pm}$. Along a continuous path $A(t)$ in $\skewm_n$, the angles can be sorted so that $\Theta(t)$ is continuous. If $\|A(0)\|_2 < l\pi$ and $\|A(1)\|_2 > 3l\pi$, then $\exists i\leq m$ such that $\theta_i(0) < l\pi$ and $\theta_i(1) > 3l\pi$. Therefore, at least one of the inequalities holds: (i) $|\theta_i(1) +\theta_j(1)| > 2l\pi$ or (ii) $|\theta_i(1) - \theta_j(1)| > 2l\pi$ for $i\neq j$. Note that $|\theta_i(0) \pm \theta_j(0)| < 2l\pi$ always holds for all $i\neq j$, which implies $A(t)$ intersects $\conjset_{l, \pm}$. In other words, every ball $\{A\in\skewm_n:\|A\|_2 > 3l\pi\}, l\in\mathbb{Z}$ separated at least two connected components, which concludes that there are infinitely many connected components in $\skewf$.
\end{proof}

In the scope of this paper, the labeling of these open components in $\skewf$ is of little interest except for the one containing the zero matrix, which is specified as the $0$-indexed component given in \cref{def:differential-skew-0-subset}.

\begin{definition}\label{def:differential-skew-0-subset}
    In particular, the $0\in\mathcal{E}$ indexed subset $\skewf_0$ is reserved for the subset that includes the zero matrix $\zerov\in\skewm_n$:
    \begin{equation}\label{eq:differential-skew-0-subset}
        \skewf_0 := \{A\in\skewm_n:\forall i\neq j, |\theta_i \pm \theta_j| < 2\pi\}.
    \end{equation}
\end{definition}

\begin{figure}[tbp]
    \centering
    \begin{subfigure}[b]{0.48\textwidth}
        \centering
        \includegraphics[width=\textwidth, height=.52\textwidth]{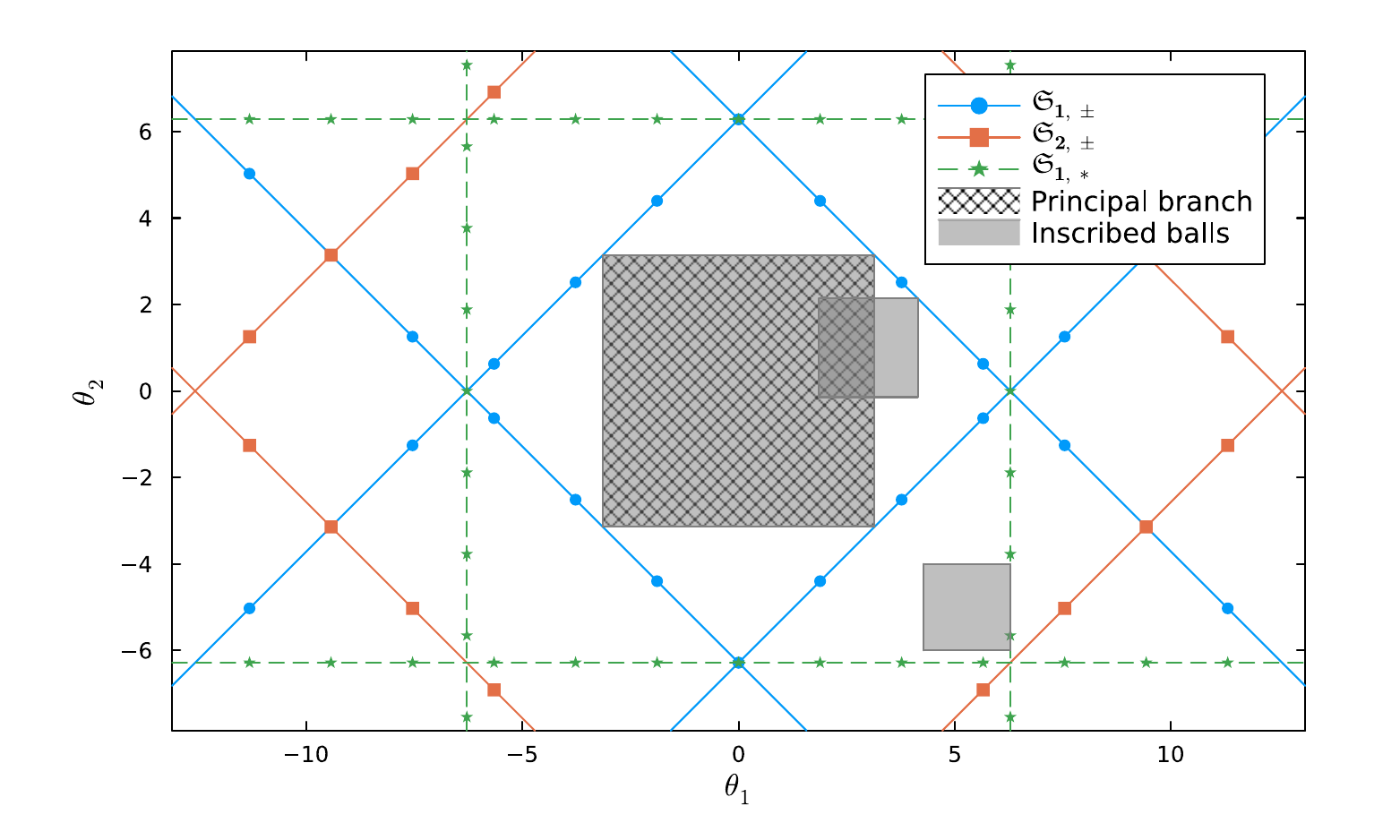}
        \caption{$(\theta_1, \theta_2)$ in $\skewm_5$.}
        \label{fig:subsets-2d}
    \end{subfigure}
    \hfill
    \begin{subfigure}[b]{0.48\textwidth}
        \centering
        \includegraphics[width=\textwidth, height=.52\textwidth]{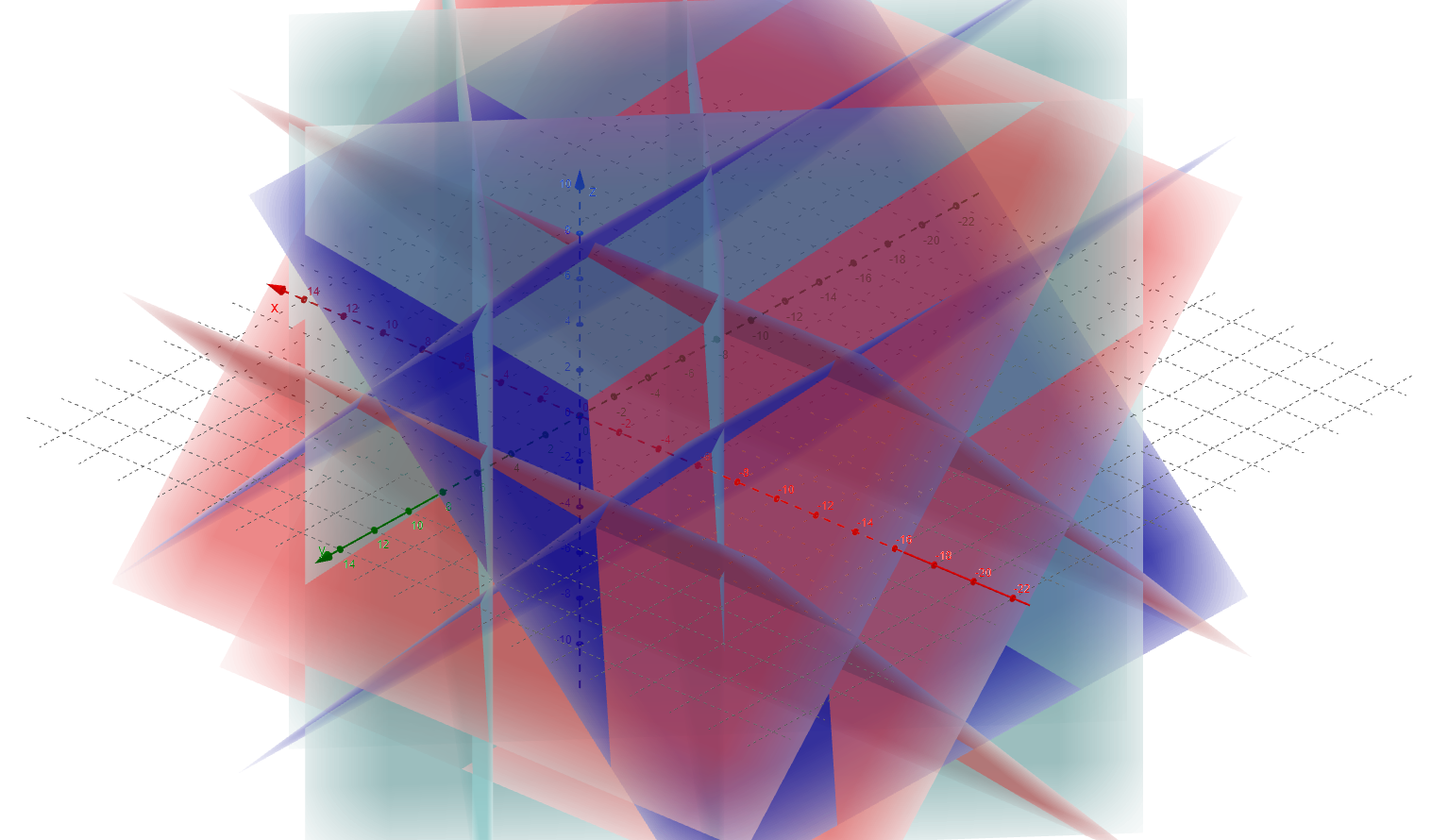}
        \caption{$(\theta_1, \theta_2, \theta_3)$ in $\skewm_6$.}
        \label{fig:subsets-3d}
    \end{subfigure}
    \caption{Illustration of $\fraks$~\cref{eq:closed-subset} and $\skewf$~\cref{eq:differential-skew-open-subsets}.}\label{fig:subsets}
\end{figure}

The image of the principal logarithm consists of matrices with eigenvalues $\lambda = a + b\im$ bounded in the imaginary part as $|b|<\pi$. For $A\in\skewm_n$ that has purely imaginary eigenvalues, the image of the principal logarithm restricted to $\skewm_n$ is
\begin{equation*}
    \ball_{\pi}(\zerov):= \{A\in\skewm_n: |\theta_i| < \pi, \forall i = 1,\ldots, m\} = \{A\in\skewm_n: \|A - \zerov\|_2 <\pi\}.
\end{equation*}
Observe that the image of the principal logarithm is contained in the $0$-indexed subset $\skewf_0$. In the $n=5$ case illustrated in \cref{fig:subsets}, $\ball_{\pi}(\zerov)$ is the largest open ball centered at $\zerov$ that \emph{touches} $\fraks$, in the sense that the ball does not intersect $\fraks$ while its closure does. In the rest of this section, the observations in the principal logarithm are generalized to a local diffeomorphism around $A\in\skewf$.

\subsection{Distance to the Tangent Conjugate Locus} Since $\fraks$ consists of closed subsets $\conjset_{l,\pm}$ and $\conjset_{l,*}$ in~\eqref{eq:closed-subset}, the distances from $A\in\skewm_n$ to these subsets and to the entire $\fraks$ are given by
\begin{equation*}
    \dist(A, \frakz) := \underset{B\in\frakz}{\min} \|A - B\|_2 \text{, for } \frakz = \conjset_{l,\pm}, \conjset_{l,*} \text{ or } \fraks.
\end{equation*}
These distances are found based on the bounds given in \cref{lemma:angular-deviation}. Consider the $2$-norm ball centered at $A$
\begin{equation*}
    \ball_{\delta}(A) := \{B\in\skewm_n: \|B-A\|_2 < \delta\}.
\end{equation*}
It follows that all $B\in\ball_{\delta}(A)$ has a vector of angles $\Sigma\in\realset^m$ that falls into the box centered at the angles $\Theta\in\realset^m$ of $A$ with side length $2\delta$, which are the shadowed boxes in \cref{fig:subsets}, including the image of the principal logarithm $ \ball_{\pi}(\zerov)$. Then, \cref{prop:distance-to-subset} finds the distance $\dist(A, \frakz)$ via the ball that touches $\frakz$.

\begin{proposition}\label{prop:distance-to-subset}
    For $A\in\skewm_n$ with Schur decomposition $(\Theta, \Xi, R)$, the distances under matrix $2$-norm from $A$ to the subsets in \cref{eq:closed-subset} are
    \begin{equation}\label{eq:distance-to-subsets}
        \lblines{
                &\dist(A,\conjset_{l,\pm}) = \underset{i\neq j}{\min} \left\{|\theta_i + \theta_j - 2l\pi|/2, |\theta_i -\theta_j -2l\pi|/2\right\},\\
                &\dist(A,\conjset_{l,*}) = \underset{i}{\min} \{|\theta_i - 2l\pi|, |\theta_i + 2l\pi|\}.\\
            }
    \end{equation}
\end{proposition}
\begin{proof}
    When $\delta = {\min}_i \{|\theta_i - 2l\pi|, |\theta_i + 2l\pi|\}$, there is $|\theta_i - 2l\pi|\geq \delta$ and $|\theta_i +2l\pi| \geq \delta$ for all $i$. It holds that $\sigma_i\neq \pm 2l\pi,\forall i$ for angles $\Sigma$ of $B\in\ball_{\delta}(A)$, i.e., $\ball_{\delta}(A) \cap \conjset_{l,*} = \emptyset$. On the other hand, let $\delta$ be obtained at the $i$th angle, i.e., one of $\theta_{i} \pm \delta$ matches $\pm 2l\pi$. Consider $B_{\pm} := \sum_{j\neq i} R_{[j]}\Xi_{[j,j]}R_{[j]}^{\T} + R_{[i]}\sbmatrix{
            0 & -\theta_{i} \mp \delta\\
            \theta_{i}\pm \delta & 0\\
        }R_{[i]}^{\T}$
    such that $\|A - B_{\pm}\|_2 = \delta$ and either $B_+ $ or $B_-$ is in $\conjset_{l,*}$, i.e., $\ball_{\delta}(A)$ is a ball touches $\conjset_{l, *}$, i.e., $\dist(A,\conjset_{l,*}) = \delta$.

    Similarly, $\delta = \underset{i\neq j}{\min} \left\{|\theta_i + \theta_j - 2l\pi|/2, |\theta_i -\theta_j -2l\pi|/2\right\}$ obtained at the $i$th and $j$th angles of $A$ also yields a ball that touches $\conjset_{l,\pm}$ at the tangency point that only have the $i$th and $j$th angles shifted by $\delta$. Note that the scalar of $1/2$ arises as the required shifting amount is evenly distributed to both angles such that one of $(\theta_{i}\pm \delta)+ (\theta_{j}\pm \delta)$ or $(\theta_{i}\pm \delta)- (\theta_{j}\mp \delta)$ matches $\pm 2l\pi$.
\end{proof}

\begin{theorem}\label{thm:distance-bound-to-conjugate-locus}
    The distance from $A\in\skewm_n$ to the tangent conjugate locus is
    \begin{equation}\label{eq:distance-to-conjugate-locus}
        \dist(A,\fraks) = \begin{cases}
            \underset{l \neq 0}{\min} \dist(A,\conjset_{l,\pm}) & n = 2m,\\
            \underset{l \neq 0}{\min} \left\{\dist(A,\conjset_{l,\pm}), \dist(A,\conjset_{l,*})\right\} & n = 2m+1,\\
        \end{cases}
    \end{equation}
    which is bounded by
    \begin{equation}\label{eq:distance-bound-to-conjugate-locus}
        \dist(A,\fraks) \leq \begin{cases}
                \pi& A\in\skewf_0,\\
                \pi/2 & A\notin \skewf_0.\\
        \end{cases}
    \end{equation}
    Hence the image of the principal logarithm $\ball_{\pi}(\zerov)$ is the largest ball that touches $\fraks$.
\end{theorem}
\begin{proof}
    Consider $A\in\skewf_0$. Since the boundary of $\skewf_0$ is contained in $\conjset_{1,\pm}$, $\dist(A,\fraks)$ must be realized by one of the following candidates
    \begin{equation*}
        \begin{matrix}
            \lblines{
                2d_{1,+} &:= \min_{i\neq j}|\theta_i+\theta_j-2\pi|\\
                2d_{1,-} &:= \min_{i\neq j}|\theta_i-\theta_j-2\pi|
            } &
            \lblines{
                2d_{-1,+} &:= \min_{i\neq j}|\theta_i+\theta_j+2\pi|\\
                2d_{-1,-} &:= \min_{i\neq j}|\theta_i-\theta_j+2\pi|
            }.
        \end{matrix}
    \end{equation*}
    If $d_{1,+} = \min\{d_{1,+}, d_{1,-}, d_{-1,+},d_{-1,-}\}$ realizes the distance, then $d_{1,+} \leq \pi$ is shown by a proof of contradiction. Suppose, contrary to the claim, that $2d_{1, +} = |\theta_{i_{\#}}+\theta_{j_{\#}}-2\pi| > 2\pi$ is realized at $i_{\#}$th and $j_{\#}$th angles. Since $|\theta_i+\theta_j| < 2\pi,\forall i, j$, then $\theta_{i_{\#}}+\theta_{j_{\#}}-2\pi < -2\pi$. It contradicts to $d_{1, +}$ being a smallest candidate as
    \begin{equation*}
        2d_{-1,+} \leq |\theta_{i_{\#}}+\theta_{j_{\#}}+2\pi| = |(\theta_{i_{\#}}+\theta_{j_{\#}} - 2\pi) + 4\pi| < 2\pi < 2d_{1,+}.
    \end{equation*}
    It is readily checked that similar contradiction is obtained for other candidates if they achieve a distance greater than $\pi$. In conclusion, $\dist(A,\conjset)\leq \pi,\forall A\in\skewf_0$. Furthermore, $\dist(A,\fraks) < \pi$ follows if $\exists |\theta_{i_{\#}}| > 0$ and $\dist(A,\fraks) = \pi$ follows if $\forall \theta_i = 0$, i.e., $A = \zerov$. It makes $\ball_{\pi}(\zerov)$ the only largest ball in $\skewf_0$.

    When $A\in\skewf_e\neq \skewf_0$, the boundary of $\skewf_e$ includes at least one of the following two pairs of subsets: (i) $\conjset_{l, \pm}$ and $\conjset_{l+1, \pm}$; (ii) $\conjset_{l, *}$ and $\conjset_{l+1, *}$, where $l, l+1\neq 0$ are integers. Apply the same technique to get $\dist(A,\skewf_e) \leq \pi/2$. Note that the larger bounds for $A\in\skewf_0$ come from the absence of the set $\conjset_{0, \pm}$ and $\conjset_{0,*}$.
\end{proof}

The ball centered at $\zerov$ that touches $\fraks$ has been shown to be the image of the principal logarithm. An important implication is that the matrix exponential, when restricted to the image of the principal logarithm, becomes a diffeomorphism. This implication can be generalized for all $A\in\skewf$ to the ball centered at $A$ that touches $\fraks$.

\subsection{Local Smooth Inversion to the Matrix Exponential}

The inverse function theorem (IFT) states that any differentiable function $f\colon \mathcal{X}\to\mathcal{Y}$ has a differentiable inverse function on $\mathcal{Y}$ if (i) $f$ is bijective and (ii) $df\colon T_x\mathcal{X} \to T_{f(x)}\mathcal{Y}$ is invertible for all $x\in\mathcal{X}$. Such a differentiable function is called a diffeomorphism.

Observe that the ball $\ball_{\pi}(\zerov)\subset \skewf_0$ is the image of the principal logarithm. In view of \eqref{eq:differential-skew-open-subsets}, $\dexpof{A}:\skewm_n\to T_{Q}\so_n$ is invertible for all $A\in \skewf$. This section further shows the bijective property of the matrix exponential in $\ball_{\pi}(A)\cap\skewf_e$ such that the matrix exponential is a local diffeomorphism, see \cref{thm:diffeomorphism}. Such a local diffeomorphism subsumes the matrix exponential on the image of the principal logarithm $\ball_{\pi}(\zerov)\cap \skewf_0 = \ball_{\pi}(\zerov)$, and the inverse matrix exponential map on a local diffeomorphism is defined as a nearby matrix logarithm, see \cref{def:nearby-lograithm}.

\begin{lemma}\label{lemma:skew-to-so-angles}
    Let $R$ be a set of Schur vectors that simultaneously block diagonalize $A = R\Xi R^{\T}\in\skewf$ and $Q = \exp(A) = RER^{\T}$. For the diagonal blocks $E_{[i,i]}, i\leq m$ and the angles $\Theta=(\theta_1,\ldots, \theta_m)\in\realset^m$, it holds that $E_{[i,i]} = E_{[j,j]} \Leftrightarrow \theta_i=\theta_j,\forall i\neq j\leq m$. Moreover, if $\exists E_{[i,i]}$ that equals $I_2$ or $-I_2$, then the following statements hold:
    \begin{enumerate}[(i)]
        \item If $n=2m+1$ and $E_{[i,i]}= I_2$, then $\theta_i = 0$;
        \item If $n=2m$ and $E_{[i,i]} = E_{[j,j]} = I_2$ for $i\neq j$, then $\theta_i=\theta_j=0$;
        \item If $n=2m$ and $\theta_i = 2l\pi$ for $l\in\mathbb{Z}\setminus \{0\}$ (consequently, $E_{[i,i]} = I_2$) then $E_{[j,j]}\neq I_2$ for $j\neq i$;
        \item If $E_{[i,i]} = -I_2$ (consequently, $\theta_i = (2l+1)\pi,l\in\mathbb{Z}$), then $E_{[j,j]}\neq -I_2$ for $j\neq i$.
    \end{enumerate}
\end{lemma}
\begin{proof}
    The proof of ``left-hand side (LHS) $\Leftarrow$ right-hand side (RHS)'' follows trivially from $E_{[i,i]} = \exp(\Xi_{i,i}),\forall i\leq m$, in which the condition $A\in\skewf$ is not used. This condition is invoked for the proof of ``LHS $\Rightarrow$ RHS''. If $E_{[i,i]} = E_{[j,j]}$, it holds that $\cos(\theta_i) = \cos(\theta_j)$ and $\sin(\theta_i) = \sin(\theta_j)$, which implies $\theta_i = \theta_j +2l\pi$ for some integer $l\in\mathbb{Z}$. In view of \cref{thm:invertibility-skew-symm}, $l = 0$ is the only candidate for $A\in\skewf$, i.e., $\theta_i = \theta_j$.

    Consider $E_{[i,i]} = I_2$ with $\theta_i = 2l_i\pi, l_i\in\mathbb{Z}$. (i) If $n = 2m+1$, then $\theta_i \neq 2l\pi, \forall 0\neq l\in\mathbb{Z}$, i.e., $\theta_i = 0$. (ii) If $n=2m$ and $E_{[j,j]} =I_2$ with $\theta_j = 2l_j\pi$, then $|\theta_i \pm\theta_j| = \|2(l_i\pm l_j)\pi\|\neq 2l\pi,\forall 0\neq l\in\mathbb{Z}$ yields $\theta_i = \theta_j = 0$. (iii) If $n = 2m$ and $\theta_i = 2l_i\pi\neq 0$, $|2l_i\pi\pm\theta_j|\neq 2l\pi,\forall 0\neq l\in\mathbb{Z}$ yields $\theta_j\neq 2l\pi,\forall l\in\mathbb{Z}$ for all $j\neq i$, i.e., $E_{[j,j]}\neq I_2, \forall i\neq j$.

    Consider $E_{[i,i]} = -I_2$ with $\theta_i = (2l_i+1)\pi,l_i\in\mathbb{Z}$. (iv) $|(2l_i+1)\pi\pm\theta_j|\neq 2l\pi,\forall 0\neq l\in\mathbb{Z}$ yields $\theta_j\neq (2l+1)\pi,\forall l\in\mathbb{Z}$ for all $j\neq i$, i.e., $E_{[j,j]}\neq -I_2, \forall i\neq j$.
\end{proof}

\begin{proposition}\label{prop:shared-schur-vectors}
    Let $(\Theta, \Xi, R_A)$ be a Schur decomposition of $A\in\skewf$, and let $A \neq B\in\skewm_n$ satisfy $\exp(A) = \exp(B) =:Q$. Then, there exists a set of Schur vectors $R_B$ of $B$ such that $R_A^{\T}QR_A = R_B^{\T}QR_B=:E$. Furthermore, $(\Theta', \Xi', R_B)$ is a Schur decomposition of $A$, i.e., $A = R_B\Xi'R_B^{\T}$, such that $\theta'_i = \theta_i + 2l_i\pi$ with $l_i\in\mathbb{Z}$ for all $i\leq m$, and that $l_i\neq 0$ only if $\theta_i = \eta\pi$ for an integer $0\neq \eta\in\mathbb{Z}$.
\end{proposition}
\begin{proof}
    Let $Q = \exp(A)$ and $E = R_A^{\T}QR_A$. Since $\exp(B) = Q$, any set of Schur vectors of $B$ is a set of Schur vectors of $Q$. For a set of Schur vectors $R_B$, one can permute the order in the columns of $R_B$ such that $R_B^{\T}QR_B = E$ and that $R_B$ remains a Schur matrix of $B$. Recall that the invariance group $\ivg(E)$ of $E$ defined in~\cref{eq:invariance-group} collects orthogonal $P$ that satisfies $P^{\T}EP = E$, then $R_A^{\T}R_B \in \ivg(E)$.

    Let $\mathcal{I}_{\mathrm{e}} = \{i \leq m: \exists \eta \in \mathbb{Z}\setminus\{0\}: \theta_i = 2\eta\pi\}$, $\mathcal{I}_{\mathrm{o}} = \{i \leq m: \exists \eta \in \mathbb{Z}: \theta_i = (2\eta+1)\pi \}$, and
    $\mathcal{I}_{\mathrm{z}} = \{i \leq m: \theta_i = 0 \}$.
    In view of \cref{lemma:skew-to-so-angles}[(iii)], either $\mathcal{I}_{\mathrm{e}}$ is an empty set or $\mathcal{I}_{\mathrm{e}}$ only has one element. If $\mathcal{I}_{\mathrm{e}}\neq \emptyset$, then $\mathcal{I}_{\mathrm{z}} = \emptyset$ as $E_{[i,i]}\neq I_2,\forall i\notin \mathcal{I}_{\mathrm{e}}$. In view of \cref{lemma:skew-to-so-angles}[(iv)], either $\mathcal{I}_{\mathrm{o}}$ is an empty set or $\mathcal{I}_{\mathrm{o}}$ only has one element.

    Observe that (i) $E_{[i,i]}$ and $\Xi_{[i,i]}$ have a different block structure (because $E_{[i,i]}$ is diagonal while $\Xi_{[i,i]}$ is not) if and only if $i\in\mathcal{I}_{\mathrm{e}}\cup\mathcal{I}_{\mathrm{o}}$; (ii) in view of \cref{lemma:skew-to-so-angles}, $\theta_i \neq \theta_j$ implies $E_{[i,i]} \neq E_{[j,j]}$; (iii) in view of \cref{thm:invertibility-skew-symm}, when $n = 2m+1$, it holds that $E_{[k,k]} = 1$, $\Xi_{[k,k]} = 0$, and, for $i<k$, $E_{[i,i]} = I_2$ if and only if $\theta_i = 0$.
    It then follows from~\cite[Lemma~A.2 and~A.3]{mataigne2024eigenvalue} and $R_A^{\T} R_B \in \ivg(E)$ that $R_A^{\T} R_B = P O$ with $P\in\ivg(\Xi )$ and $O = \mathrm{diag}(O_{[1,1]},\dots,O_{[k,k]})$, where $O_{[i,i]}$ is a $2\times 2$ rotation or reflection matrix if $i\in\mathcal{I}_{\mathrm{e}}\cup\mathcal{I}_{\mathrm{o}}$ and $O_{[i,i]}$ is the identity matrix of the same size (2, or 1 if $i = k > m$) otherwise. Hence $R_B^{\T}A R_B = R_B^{\T}R_A \Xi R_A^{\T}R_B = (PO)^{\T}\Xi PO = O^{\T}\Xi O = \Xi'$ where $\Xi'_{[i,i]} = -\Xi_{[i,i]}$ if $O_{[i,i]}$ is a reflection and $\Xi'_{[i,i]} = \Xi_{[i,i]}$ otherwise. Hence $\theta'_i = \theta_i$ or $\theta'_i = -\theta_i \neq 0$, and the latter happens only if $O_{[i,i]}$ is a reflection, which happens only if $i\in\mathcal{I}_{\mathrm{e}}\cup\mathcal{I}_{\mathrm{o}}$, i.e., $\theta_i = \eta\pi$ with $\eta\in\mathbb{Z}_0$, yielding $\theta'_i = \theta_i - 2\eta\pi$.
\end{proof}

\begin{proposition}\label{prop:isolated-preimage}
    If $A\in\skewf$ and $A\neq B\in\skewm_n$ that satisfy $\exp(A) = \exp(B)$, then $A$ and $B$ are separated by a distance of at least $2\pi$ in matrix $2$-norm, i.e.,
    \begin{equation}\label{eq:isolated-preimage-distance}
        \|A-B\|_2 \geq 2\pi,\forall A \in\skewf, A\neq B\in\skewm_n, \exp(A) = \exp(B).
    \end{equation}
\end{proposition}
\begin{proof}
        Let $(\Theta, \Xi, R_A)$ and $(\Sigma, H, R_B)$ be the Schur decomposition of $A$ and $B$ such that $E := \exp(\Xi ) = \exp(H)$. Since $\cos(\theta_i) = \cos(\sigma_i),\sin(\theta_i) = \sin(\sigma_i), \forall i \leq m$, each angle is differed by a multiple of $2\pi$, i.e., $\exists l_i \in\mathbb{Z}$ such that $\theta_i - \sigma_i=2l_i\pi,\forall i\leq m$. In view of \cref{prop:shared-schur-vectors}, $(\Theta', \Xi', R_B)$ is a Schur decomposition of $A$ such that $\theta_i' = \theta_i + 2l_i'\pi$ with $l'_i \in\mathbb{Z}$ for all $i\leq m$.

       Notice that $\|A-B\|_2 = \|R_B^{\T}(A-B)R_B\|_2 = \|\Xi' - H\|_2 = \max_{i\leq m}|\theta_i' - \sigma_i|$, where $|\theta_i'-\sigma_i| = |2l_i\pi + 2l_i'\pi|$ holds for all $i\leq m$. Therefore, $\|A-B\|_2 = \max_{i\leq m}|\theta'_i - \sigma_i|$ is a multiple of $2\pi$. Since $A\neq B \Leftrightarrow \|A-B\|_2 \neq 0$, it holds that $\|A-B\|_2 \geq 2\pi$.
\end{proof}

The separation by at least $2\pi$ implies the bijection on $\exp:\skewf_e\cap\ball_{\pi}(A)\to \so_n,\forall A\in\skewf_e$, which leads to a local well-defined inversion given in \cref{thm:diffeomorphism}.

\begin{theorem}\label{thm:diffeomorphism}
    For $A\in\skewf_e\subset \skewf$, the matrix exponential
    \begin{equation*}
        \exp: \calm_A:=\skewf_e\cap \ball_{\pi}(A) \to\{\exp(X):X\in\calm_A\}=: \calz_A
    \end{equation*}
    is a diffeomorphism that has a smooth inversion. Moreover, $\ball_{\dist(A,\fraks)}(A)\subset \calm_A$.
\end{theorem}
\begin{proof}
    For all $X,Y\in\ball_{\pi}(A)$, the triangle inequality shows that
    \begin{equation*}
        \|X-Y\|_2 \leq \|X - A\|_2+\|Y-A\|_2 < 2\pi,\forall X,Y\in\ball_{\pi}(A).
    \end{equation*}
    According to \cref{prop:isolated-preimage}, $\exp(X) \neq \exp(Y)$ holds $\forall A, B\in \calm_A\subset \skewf_e$. In other words, $\exp:\calm_A\to\calz_A$ is a smooth bijection with invertible differentiation over the entire $\calm_A$. Then, IFT guarantees that the inversion is also a smooth bijection with invertible differentiation over the entire $\calz_A$. Finally, $\ball_{\dist(A,\fraks)}(A)\subset \calm_A$ follows from $\ball_{\dist(A,\fraks)}(A)\subset \skewf_e$ and $\dist(A,\fraks) \leq \pi$.
\end{proof}

\begin{definition}\label{def:nearby-lograithm}
    The inverse of the smooth bijection $\exp:\calm_A\to\calz_A$ is termed the \emph{nearby matrix logarithm} around $A$, denoted by
    \begin{equation}\label{eq:nearby-logarithm}
        \nlog_A:\calz_A\to \calm_A.
    \end{equation}
\end{definition}

As pointed out earlier in this section, the exponential on $\ball_{\pi}(\zerov) = \ball_{\pi}(\zerov)\cap \skewf_0$ is shown to be a local diffeomorphism. In view of \cref{def:nearby-lograithm}, the principal logarithm is the nearby matrix logarithm around $\zerov$. Note that the matrix exponential on the entire $\skewf_0$ is not a diffeomorphism due to the lack of the bijection on it.

\subsection{Continuous Motions on Skew-Symmetric Matrices}

As an important theoretical application of the invertibility of $\dexpof{A}:\skewm_n\to T_Q\so_n$, the nearby matrix logarithm~\cref{eq:nearby-logarithm} can be used to identify a continuous motion $A(t)$ in $\skewm_n$ that captures $Q(t) = \exp(A(t))$ given in $\so_n$. The detailed discussion on computing the nearby matrix logarithm is beyond the scope of this paper. For the purposes of demonstrating theoretical significance in $\nlog$, an algorithm for computing $\nlog$ proposed in~\cite[Algorithm 4]{deng2024smoothly} is employed to compute $A(t)$ from some simple $Q(t)$.

There are three smooth curves on $\so_4$ investigated in this experiment, all in the form of $Q(t) = Q \exp(t Y), t\in [0, 1], Q\in\so_n, Y\in\skewm_n$.
\begin{enumerate}
    \item $Q_1(t)$ is generated by random $Q_1$ and $Y_1$ with $A_1(0) = \log(Q_1)$.
    \item $Q_2(t)$ is generated by random $Q_2$ and $Y_2 = R\diagm\left(\zerov, \Xi_{[2,2]}\neq \zerov\right)R^{\T}$ with $A_2(0) = \log(Q_2)$, where $R$ is a set of Schur vectors of $Q_2$.
    \item $Q_3(t)$ is chosen such that $A_3(t_*)\in\fraks$ with $\dexpof{A_3(t_*)}[X_3] = Q_3(t_*)Y_3$.
\end{enumerate}
The curve $Q_2(t)$ is designed to have a simple motion while the curve $Q_3(t)$ is a designed corner case that moves across the $\{\exp(B):B\in\fraks\}$. More specifically, $t_*\in[0,1]$ and $A_3(t_*)\in\fraks, X_3\in\skewm_n$ are chosen first. Then, $Q_3(t_*) = \exp(A_3(t_*))$ and $Q_3(t_*)Y_3 = \dexpof{A_3(t_*)}[X_3]$ are computed. Finally, $Q_3(0)$ are computed long the reverse direction: $Q_3'(t) = Q_3(t_*)\exp(-t Y_3), t\in [0, t_*]$, such that $Q_3(0) = Q_3'(t_*)$.

\begin{figure}[tbp]
    \centering
    \includegraphics[width=.6\textwidth]{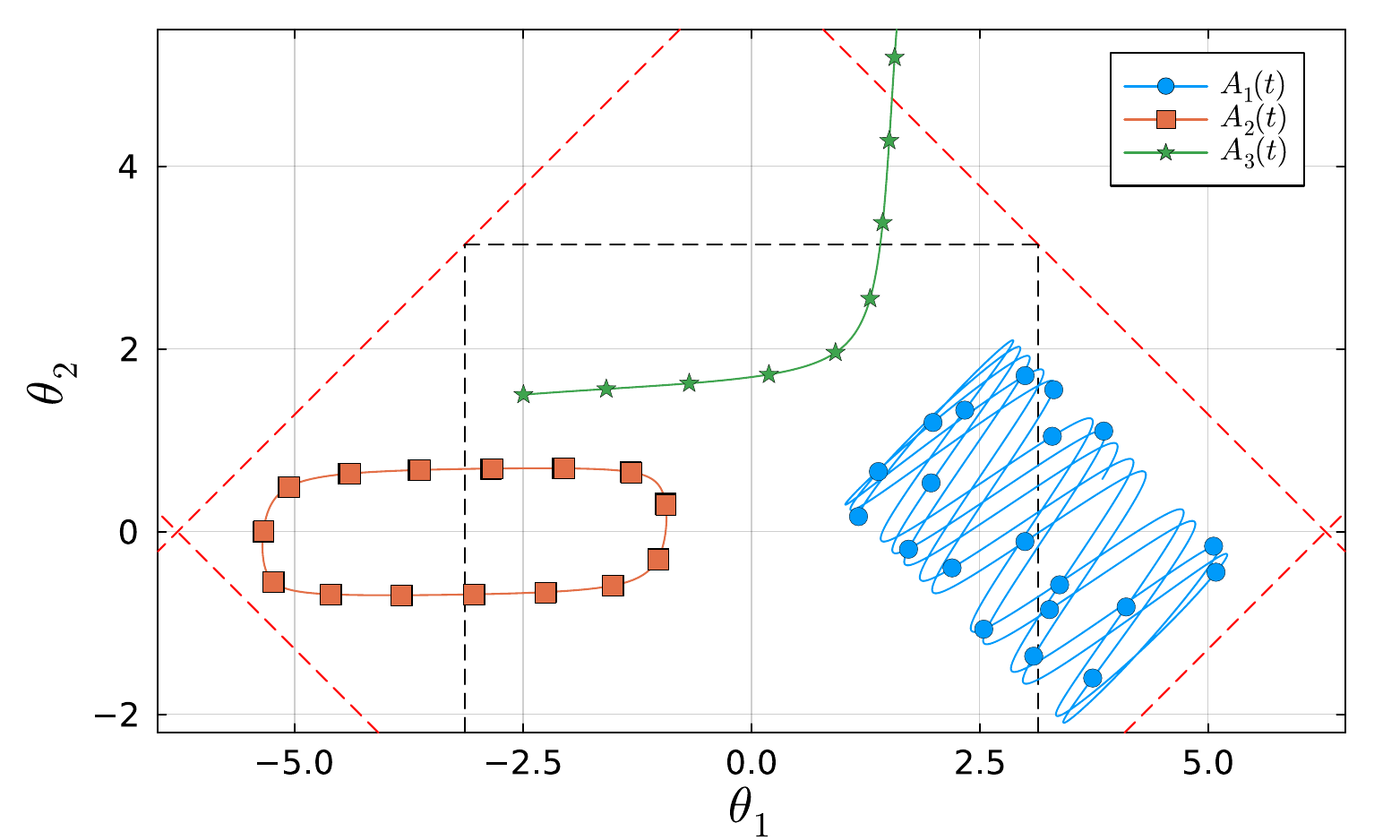}
    \caption{Curves $A(t)$ in $\skewm_4$ that satisfy $\exp(A(t)) = Q \exp(t Y)$ with $Y\in\skewm_4$.}\label{fig:st-curve}
\end{figure}

Since the matrices in $\skewm_n$ are difficult to visualize, the angular plot in \cref{fig:subsets} is used for visualization. A Schur decomposition $(\Theta(t), D(t), R(t))$ is continuously assigned to the computed $A(t)\in \skewm_n$ such that $\Theta(t) = (\theta_1(t),\theta_2(t))$ reported in \cref{fig:st-curve} form a continuous curve for each instance.

The curves $A_1(t)$ and $A_2(t)$ both move within $\skewf_0$ but beyond the image of the principal logarithm $\ball_{\pi}(\zerov)$. Note that $A_1(t)$ does not intersect with itself as $R_1(t)$ is also varying along $t$ but $A_2(t)$ forms a closed loop as $R_2(t)$ remains constant along $t$. The corner case $A_3(t)$ has crossed the $\fraks$ numerically as expected. The ``numerically crossing'' here means that there is a tiny interval $t_{i}\lessapprox t_* \lessapprox t_{i+1}$ where $A_3(t_i)\in \skewf_0$ and $A_3(t_{i+1})\notin\skewf_{\zerov}$. Then $A_3(t)$ crosses the boundary $\fraks$, presumably at $t = t_*$.

\section{Numerical Experiments on the Differentiations}\label{sec:experiments}

This section reports experiments on the formulae of $\dexpof{A}[X] = Q Y$ and $\dexpinvof{A}[Q Y] = X$ for $A,X,Y\in\skewm_n$. The exponential $Q=\exp(A)$ is also included in the experiments as an important preprocessing task that affects the overall compute time.  In particular, the implementations being tested include:
\begin{itemize}
    \item $\dexpof{A}[X] = Q Y$: The proposed formula~\eqref{eq:dskew-skew}; The Pad{\'e} approximant proposed in~\cite{al2009exponential}; Daletski$\breve{\mathrm{\i}}$--Kre$\breve{\mathrm{\i}}$n formulae~\eqref{eq:dexp-diag} and~\eqref{eq:linear-diag}.
    \begin{itemize}
        \item Skew-Symmetric Output $X\mapsto Y$.
        \item Full Differentiation Output $X\mapsto Q Y$.
    \end{itemize}
    \item $\dexpinvof{A}[Q Y] = X$: The proposed formula~\eqref{eq:dexp-skew-inv}; The Pad{\'e} approximant proposed in~\cite{al2013logarithm}; Inversion of Daletski$\breve{\mathrm{\i}}$--Kre$\breve{\mathrm{\i}}$n formulae~\eqref{eq:dexp-diag} and~\eqref{eq:linear-diag}.
    \begin{itemize}
        \item Skew-Symmetric Input $Y\mapsto X$.
        \item Full Differentiation Input $Q Y\mapsto X$.
    \end{itemize}
    \item $\exp(A) = Q$ computed by the following formulae:
    \begin{itemize}
        \item Schur decomposition $Q = R\exp(\Xi )R^{\T}$,
        \item Pad{\'e} Approximant $\poly{r/r}{\exp}(A)$ with designated order,
        \item \texttt{MATLAB} Built-in Function \texttt{expm} (Adaptive Pad{\'e} Approximant).
    \end{itemize}
\end{itemize}
The order $[r/r]$ Pad{\'e} approximants are denoted by ``Pade-[r/r]''. The proposed formula for skew-symmetric matrices are denoted as ``Skew-Symm'' and Daletski$\breve{\mathrm{\i}}$--Kre$\breve{\mathrm{\i}}$n formulae for semi-simple matrices are denoted as ``Semi-Simple''.

\begin{figure}[tbp]
    \centering
    \begin{subfigure}[b]{0.48\textwidth}
        \centering
        \includegraphics[width=\textwidth, height=.52\textwidth]{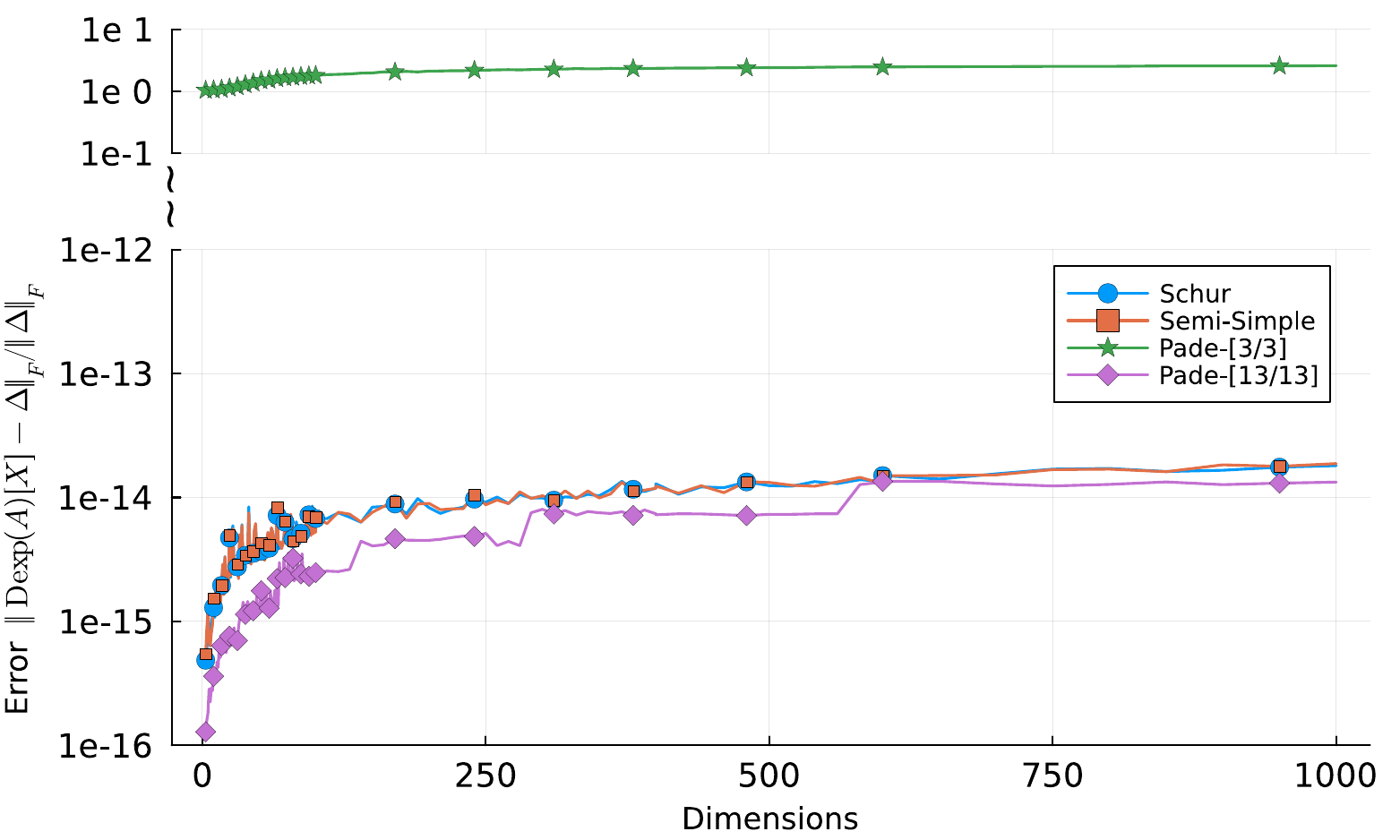}
        \caption{Error of $X\mapsto \dexpof{A}[X]$}
        \label{fig:forward-error}
    \end{subfigure}
    \hfill
    \begin{subfigure}[b]{0.48\textwidth}
        \centering
        \includegraphics[width=\textwidth, height=.52\textwidth]{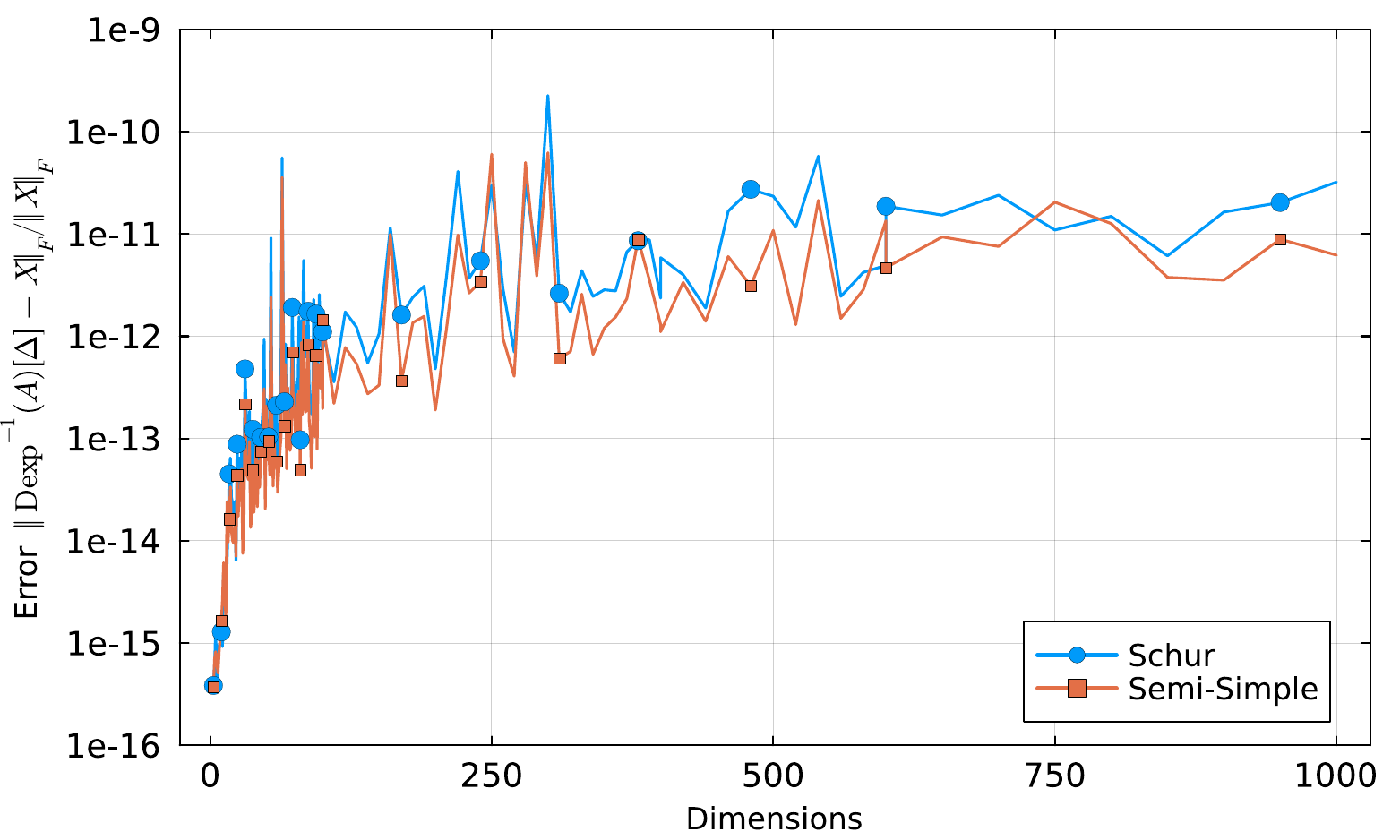}
        \caption{Error of $\Delta\mapsto \dexpinvof{A}[\Delta]$}
        \label{fig:inverse-error}
    \end{subfigure}
    \caption{Errors of $\dexpof{A}[X]=\Delta$ (Compared with $256$-bits Floating Point System).}
    \label{fig:error}
\end{figure}

Every skew-symmetric input is generated with random numbers uniformly sampled on $[-1, 1]$. For every routine timed, one warm-up run is performed before $20$ or more executions. The elapsed consumptions of these executions are recorded and the average is reported as the final timed elapsed consumption. For the error reported in \cref{fig:error} and \cref{fig:expm-error}, the numerically accurate matrices in $\dexpof{A}[X]=\Delta$ and $\exp(A) = Q$ are computed using the $256$-bit floating point system natively supported in \texttt{Julia} via \texttt{BigFloat} that keeps $70$ decimal digits of precision, while other computations are performed in the standard double precision floating point system ($64$ bits, about $15$ decimal digits of precision). In particular, each double-precision matrix is converted to \texttt{BigFloat} (by attaching \texttt{0}), computed in the $256$-bit arithmetic, and then cast back to double precision (by truncating tails) to obtain the numerically ``accurate'' result. The error is measured by $\mathrm{error}(M) = \|M - M_*\|_{\F}/\|M_*\|_{\F}$ where $\|\cdot\|_{\F}$ denotes the Frobenius norm and $M_*$ denotes the accurate matrix. Note that the accurate $\dexpof{A}[X]$ uses the divided difference $\exp(A+\varepsilon X) - \exp(A) / \varepsilon$ with $\varepsilon = \texttt{1e-20}$ and the exponential uses the $[13/13]$ order Pad{\'e} approximant. The source code used in this section is publicly available at \url{https://github.com/zhifeng1703/cpp-released-code/dexp}. The utility routines, including the Schur decomposition, the matrix exponential, and the differentiation formulae studied in this paper, are compiled into a \texttt{MATLAB}-callable library together with a ready-to-use wrapper, all of which are provided in the same repository.

\subsection{Computing Environment}
All experiments were conducted on \texttt{Windows Subsystem Linux} and executed in \texttt{C++} (except for the \texttt{MATLAB} (R2024a) built-in \texttt{expm}).
{\footnotesize
\begin{itemize}
    \item CPU/RAM: \texttt{13th Gen Intel(R) Core(TM) i5-13600KF, 14 Cores, 3.50 GHz} / \texttt{16 GB}.
    \item OS: \texttt{Windows 11, Windows Subsystem Linux (Ubuntu 22.04.5 LTS)}.
    \item BLAS and LAPACK Library: \texttt{Intel(R) Math Kernel Library}.
\end{itemize}
}

\subsection{Numerical Results}\label{subsec:results}

Overall, the proposed formulae and their complex-origin counterparts, the Daletski$\breve{\mathrm{\i}}$--Kre$\breve{\mathrm{\i}}$n formulae, are accurate up to sufficient number of decimal digits ($13\text{--} 15$ digits in $X\mapsto \dexpof{A}[X]$ and $10\text{--} 15$ digits in $\Delta\mapsto \dexpinvof{A}[\Delta]$). The order $[13/13]$ Pad{\'e} approximant obtains slightly better precision in $X\mapsto \dexpof{A}[X]$, while the order $[3/3]$ Pad{\'e} approximant obtains almost no meaningful digit, as it is designed for sufficiently small $\|A\| < \texttt{1.08e-2}$ without squaring, \cite[Table 6.1]{al2009exponential}. Note that the Pad{\'e} approximant of $\der\log(\exp(A))[\Delta]$ does not return $X$ when $A$ is not in the image of the principal logarithm. In particular, the randomly sampled $A\in\skewm_n$ with entries uniformly distributed in $[-1,1]$ almost always lies outside the image of the principal logarithm when $n > 20$. Therefore, the corresponding errors are not reported in \cref{sub@fig:inverse-error}. In terms of computational complexity, we observe the characteristic fluctuations in elapsed time caused by the parity of $n$, since odd values of $n$ require additional $O(n^2)$ ops. Note that the proposed formulae are significantly faster than the existing robust formulae. While the order $[3/3]$ Pad{\'e} approximant of $\dexpof{A}$ is faster in one call to the full differentiation, it is applicable only for sufficiently small $A$, as shown in \cref{sub@fig:forward-error}. Moreover, the advantage by about $0.1$ ms across all dimensions becomes less desirable as the relative speedup decreases from nearly twice as fast for $n < 40$ to almost no difference for $n > 250$. With $n > 250$, the numerical experiments support the complexity analyses in \cref{tab:dexp-ops} and \cref{tab:dlog-ops}, as the speedups in the recorded time are consistent with the multiples of the $n^3$ ops. As the problem size becomes smaller, the $O(n^2)$ terms in the complexity analyses become significant. Therefore, the $O(n^3)$ term no longer accurately predicts the trends of compute time reported in \cref{fig:timing_dexp_skew}$\sim$\cref{fig:timing_dlog_full} with $n < 50$.

\Cref{fig:timing_dexp_skew}$\sim$\Cref{fig:timing_dlog_full} report two different elapsed times: the subgraph on the left reports the total elapsed time of one execution and the subgraph on the right reports the elapsed time broken down into the preprocessing stage and the computing stage, where the preprocessing stage can be reused in multiple calls. The computing stage is reported above the $x$-axis and the preprocessing stage is reported below the $x$-axis, which is growing downward. The region in between is exactly the total elapsed time reported on the left, i.e., the graph on the right is shifting the graph on the left by the amount of elapsed time consumed in the preprocessing stage. Consider $\dexpof{A}[X]$ reported in \cref{fig:timing_dexp_skew} and \cref{fig:timing_dexp_full} as an example. When there is only one call to $\dexpof{A}[X]$, the graph on the left is sufficient. When there are multiple calls to $\dexpof{A}[X]$ for the same $A$ with different $X$'s, the computing stage on the right becomes more dominant in the overall complexity as more calls are made.

\subsubsection{Differentiation}

\begin{figure}[tbp]
    \centering
    \begin{subfigure}[b]{0.48\textwidth}
        \centering
        \includegraphics[width=\textwidth, height=.52\textwidth]{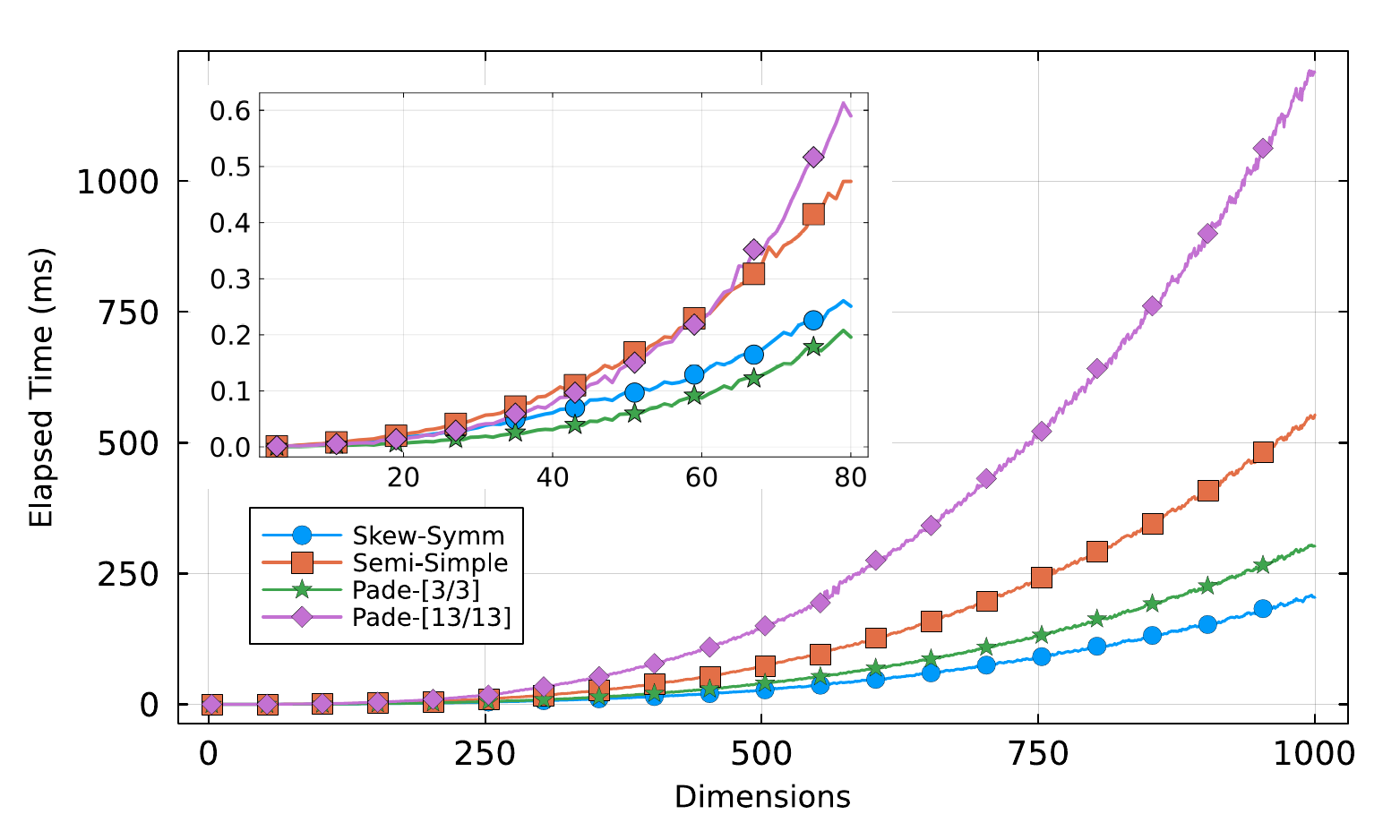}
        \caption{Total time}
        \label{fig:total_time_dexp_skew}
    \end{subfigure}
    \hfill
    \begin{subfigure}[b]{0.48\textwidth}
        \centering
        \includegraphics[width=\textwidth, height=.52\textwidth]{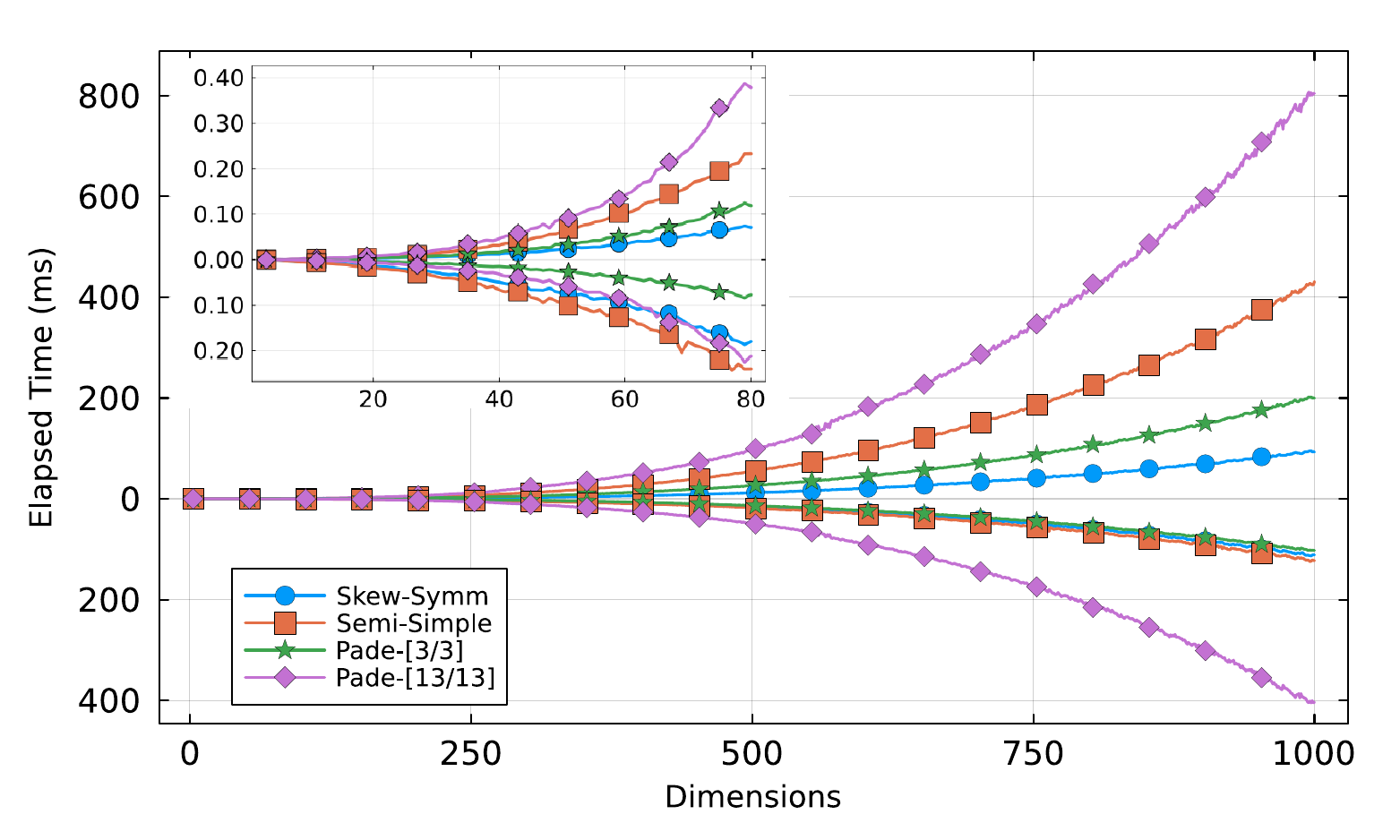}
        \caption{Breakdown (Preprocessing below $x$-axis)}
        \label{fig:breakdown_dexp_skew}
    \end{subfigure}
    \caption{Execution time for $X\mapsto Y$ in $\dexpof{A}[X] = Q Y$ across algorithms.}\label{fig:timing_dexp_skew}
\end{figure}

\begin{figure}[htbp]
    \centering
    \begin{subfigure}[b]{0.48\textwidth}
        \centering
        \includegraphics[width=\textwidth, height=.52\textwidth]{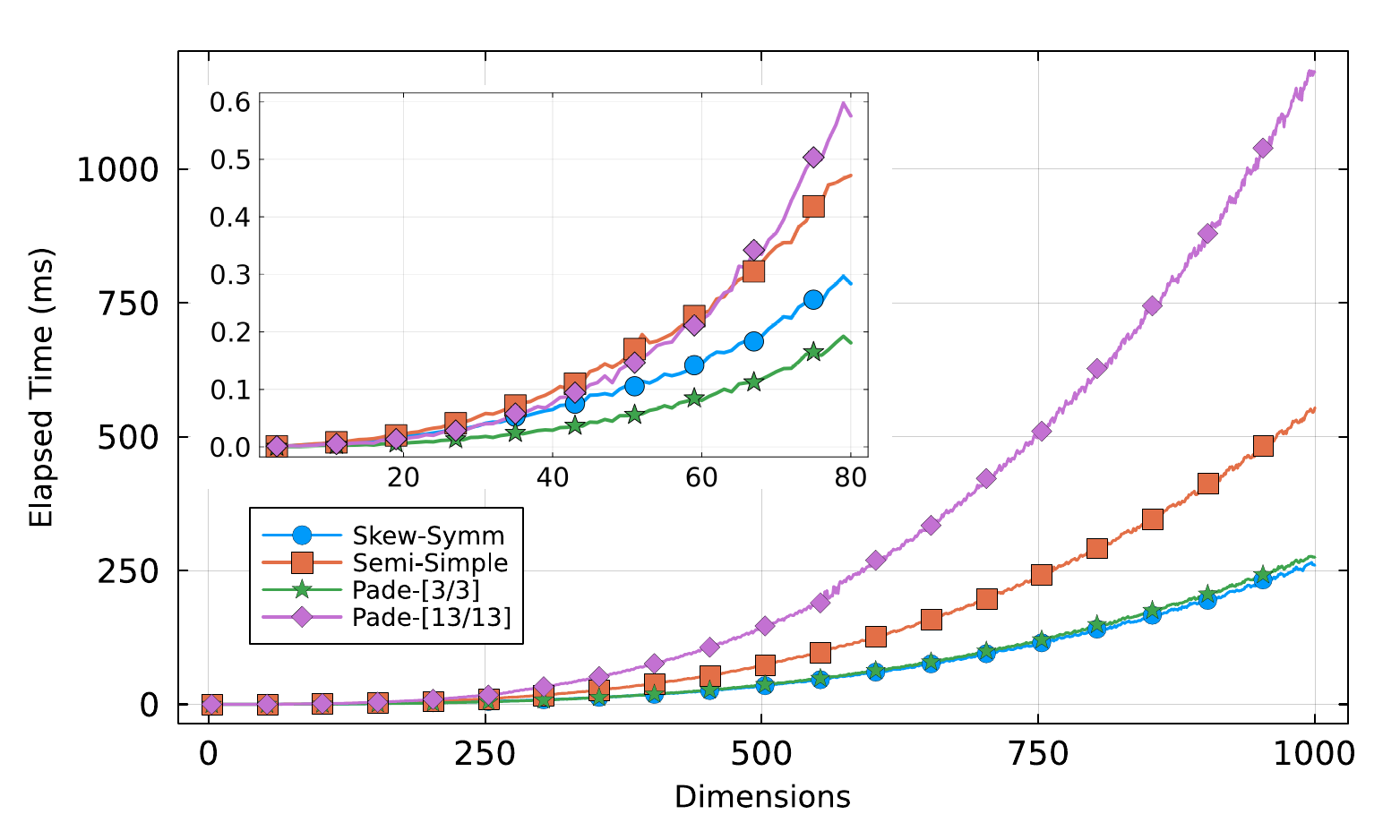}
        \caption{Total time}
        \label{fig:total_time_dexp_full}
    \end{subfigure}
    \hfill
    \begin{subfigure}[b]{0.48\textwidth}
        \centering
        \includegraphics[width=\textwidth, height=.52\textwidth]{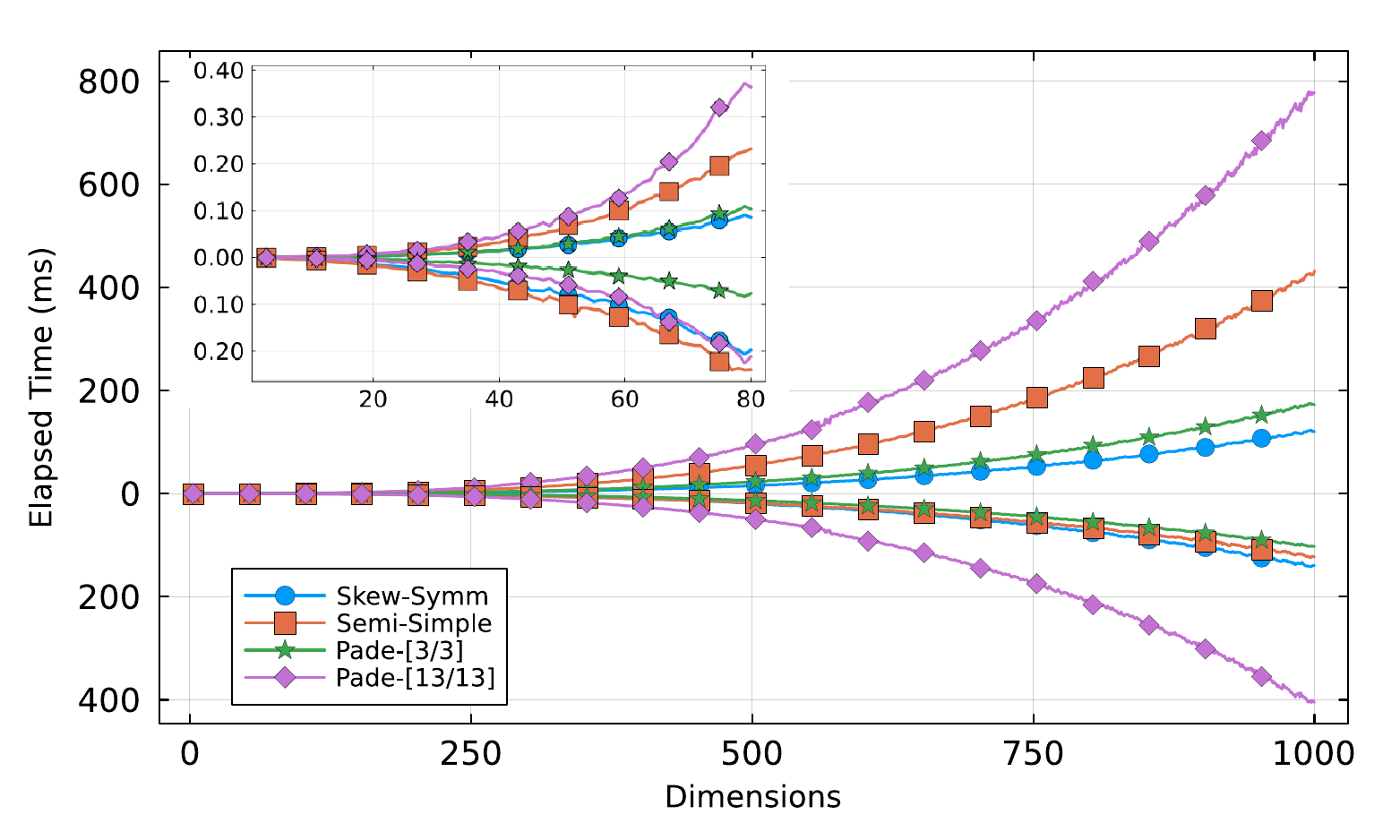}
        \caption{Breakdown (Preprocessing below $x$-axis)}
        \label{fig:breakdown_dexp_full}
    \end{subfigure}
    \caption{Execution time for $X\mapsto \dexpof{A}[X]$ across algorithms.}\label{fig:timing_dexp_full}
\end{figure}

\cref{fig:timing_dexp_skew} and \cref{fig:timing_dexp_full} present the efficiency comparison in computing $X\mapsto Y$ and $X\mapsto Q Y$ from $\dexpof{A}[X]=Q Y$. Note that the computations of input parameters like the Schur/eigen-decomposition of $A$ and the exponential $Q=\exp(A)$ are included in the comparison. An excerpt of the recorded elapsed times is given in \cref{tab:dexp-elapsed}.

Consider $n > 250$. For one execution of $X\mapsto \dskew_A(X)$, reported in \cref{sub@fig:total_time_dexp_skew}, the proposed formula is $2.3$ times faster than Daletski$\breve{\mathrm{\i}}$--Kre$\breve{\mathrm{\i}}$n formula. For the dominant computing stage in multiple calls of $X\mapsto \dskew_A(X)$, reported in \cref{sub@fig:breakdown_dexp_skew}, the proposed formula is $3.9$ times faster than Daletski$\breve{\mathrm{\i}}$--Kre$\breve{\mathrm{\i}}$n formula. For one execution of $X\mapsto \dexpof{A}[X]$, reported in \cref{sub@fig:total_time_dexp_full}, the proposed formula is $2.2$ times faster than Daletski$\breve{\mathrm{\i}}$--Kre$\breve{\mathrm{\i}}$n formula. For the dominant computing stage in multiple calls of $X\mapsto \dexpof{A}[X]$, reported in \cref{sub@fig:breakdown_dexp_full}, the proposed formula is $3.7$ times faster than Daletski$\breve{\mathrm{\i}}$--Kre$\breve{\mathrm{\i}}$n formula.

\begin{table}[tbp]
    \centering
    \caption{Recorded elapsed time for $Y = \dskew_A(X)$ and $\Delta=\dexpof{A}[X]$}\label{tab:dexp-elapsed}
    \tiny{
        \begin{tabular}{|l|ccc|ccc|ccc|}
            \hline
            $X\mapsto Y$& \multicolumn{3}{c|}{Preprocessing} & \multicolumn{3}{c|}{Computation} & \multicolumn{3}{c|}{Total} \\
            Size $n$& $10$ & $50$ & $1000$& $10$ & $50$ & $1000$& $10$ & $50$ & $1000$\\
            \hline
            Skew-Symm& 5.3e-3& 7.5e-2& 1.1e+2& \textbf{1.2e-3}& \textbf{2.3e-2}& \textbf{9.3e+1}& 6.5e-3& 9.8e-2& \textbf{2.0e+2}\\
            Pad{\'e}-$[3/3]$& \textbf{1.3e-3}& \textbf{2.6e-2}& \textbf{1.0e+2}& 1.5e-3& 3.3e-2& 2.0e+2& \textbf{2.7e-3}& \textbf{5.9e-2}& 3.0e+2\\
            Pad{\'e}-$[13/13]$& 2.1e-3& 5.5e-2& 4.0e+2& 2.6e-3& 9.1e-2& 8.0e+2& 4.7e-3& 1.5e-1& 1.2e+3\\
            Semi-Simple& 6.0e-3& 9.7e-2& 1.2e+2& 1.3e-3& 6.2e-2& 4.3e+2& 7.3e-3& 1.6e-1& 5.5e+2\\
            \hline
            $X\mapsto \Delta$& \multicolumn{3}{c|}{Preprocessing} & \multicolumn{3}{c|}{Computation} & \multicolumn{3}{c|}{Total} \\
            Size $n$& $10$ & $50$ & $1000$& $10$ & $50$ & $1000$& $10$ & $50$ & $1000$\\
            \hline
            Skew-Symm& 5.6e-3& 8.0e-2& 1.4e+2& 1.3e-3& \textbf{2.7e-2}& \textbf{1.2e+2}& 6.8e-3& 1.1e-1& \textbf{2.6e+2}\\
            Pad{\'e}-$[3/3]$& \textbf{1.3e-3}& \textbf{2.6e-2}& \textbf{1.0e+2}& 1.4e-3& 2.8e-2& 1.7e+2& \textbf{2.6e-3}& \textbf{5.4e-2}& 2.7e+2\\
            Pad{\'e}-$[13/13]$& 2.1e-3& 5.5e-2& 4.0e+2& 2.5e-3& 8.7e-2& 7.8e+2& 4.6e-3& 1.4e-1& 1.2e+3\\
            Semi-Simple& 6.0e-3& 9.5e-2& 1.2e+2& \textbf{1.2e-3}& 6.1e-2& 4.3e+2& 7.2e-3& 1.6e-1& 5.5e+2\\
            \hline
            \multicolumn{10}{l}{Skew-Symm: Formula~\eqref{eq:dskew-skew} \:\:  Pad{\'e}-$[r/r]$: Pad{\'e} approximant~\cite{al2009exponential} \:\: Semi-Simple: Formula~\eqref{eq:dexp-diag} and~\eqref{eq:linear-diag}}
        \end{tabular}
    }
\end{table}

\subsubsection{Inverse Differentiation}

\begin{figure}[tbp]
    \centering
    \begin{subfigure}[b]{0.48\textwidth}
        \centering
        \includegraphics[width=\textwidth, height=.52\textwidth]{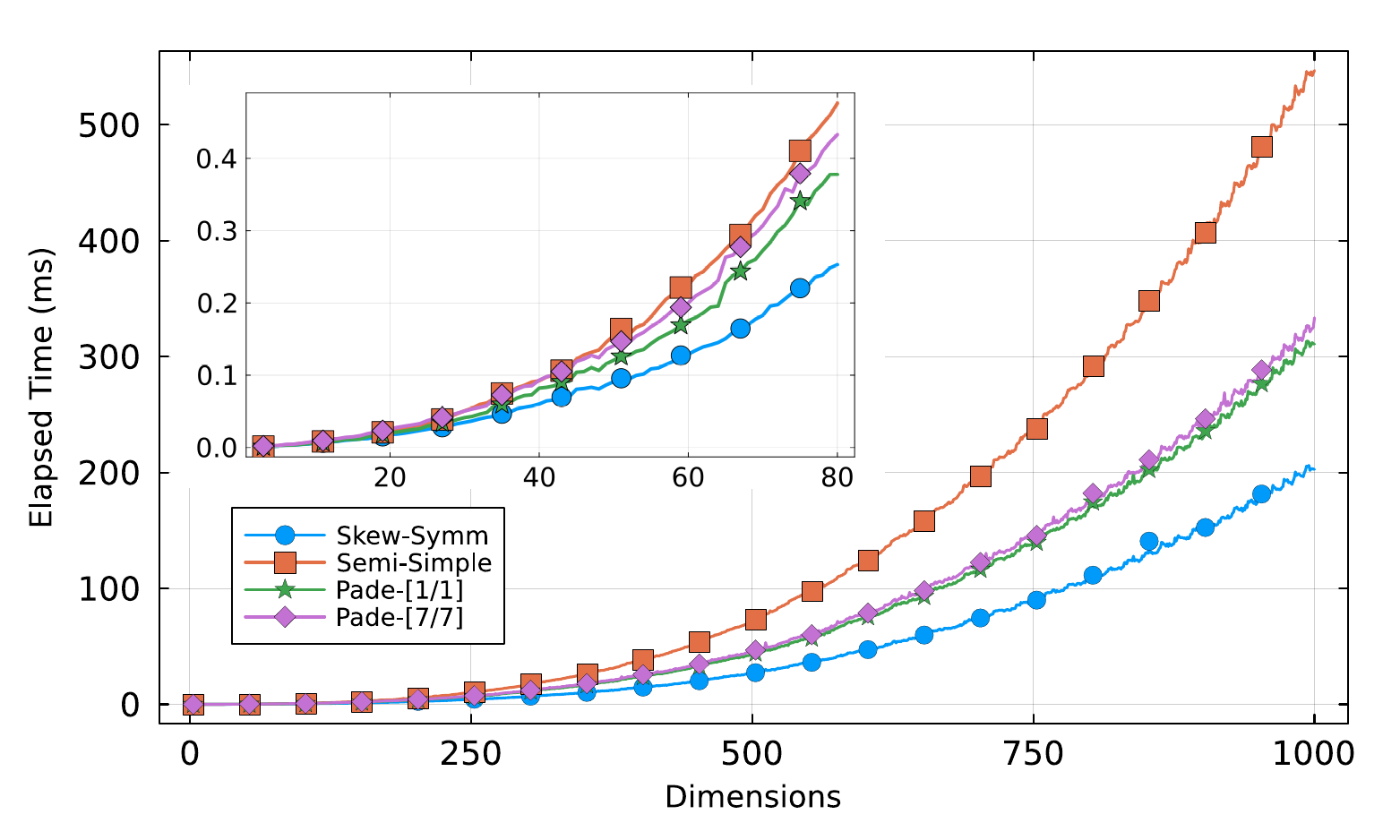}
        \caption{Total time}
        \label{fig:total_time_dlog_skew}
    \end{subfigure}
    \hfill
    \begin{subfigure}[b]{0.48\textwidth}
        \centering
        \includegraphics[width=\textwidth, height=.52\textwidth]{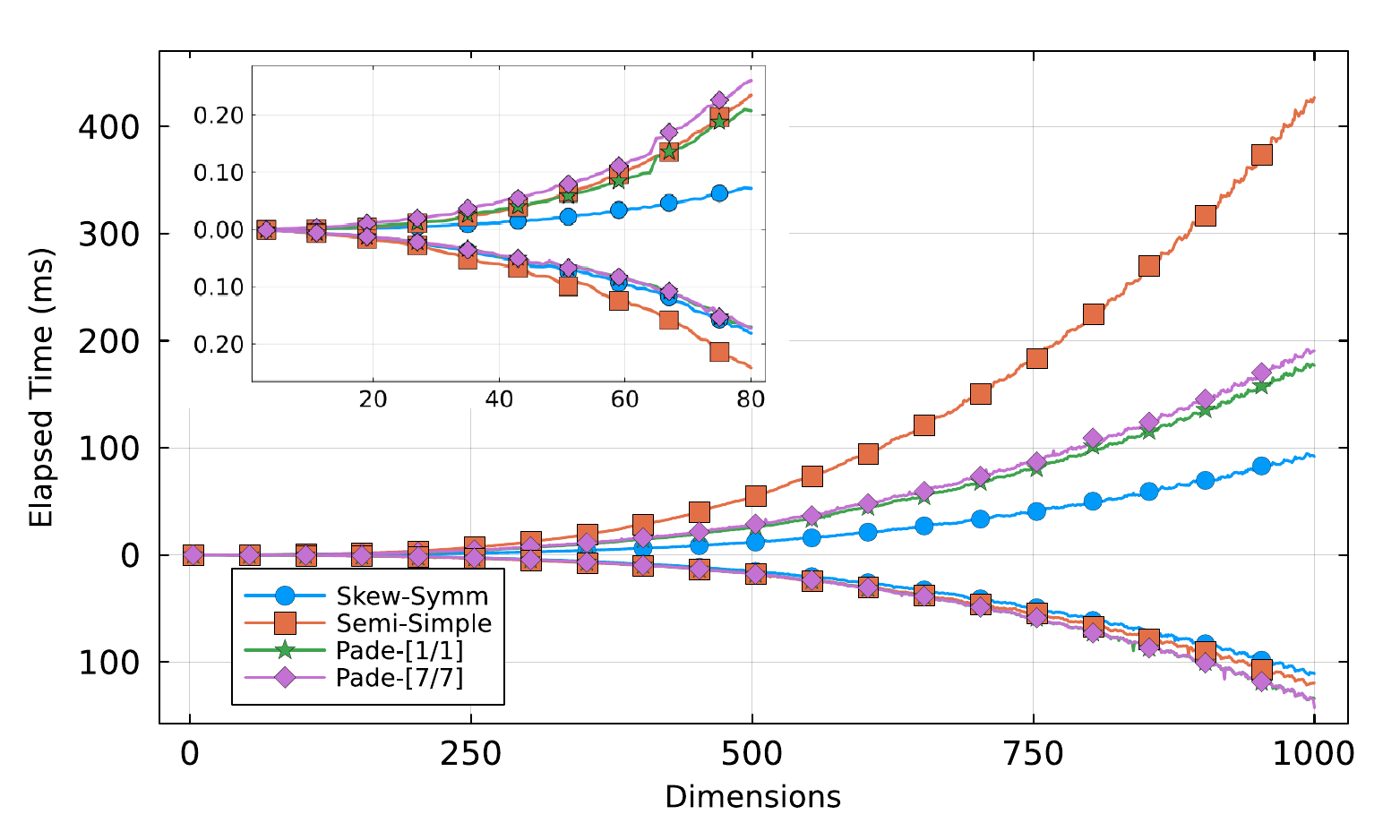}
        \caption{Breakdown (Preprocessing below $x$-axis)}
        \label{fig:breakdown_dlog_skew}
    \end{subfigure}
    \caption{Execution time for $Y\mapsto X$ in $\dexpinvof{A}[Q Y] = X$ across algorithms.}\label{fig:timing_dlog_skew}
\end{figure}

\begin{figure}[tbp]
    \centering
    \begin{subfigure}[b]{0.48\textwidth}
        \centering
        \includegraphics[width=\textwidth, height=.52\textwidth]{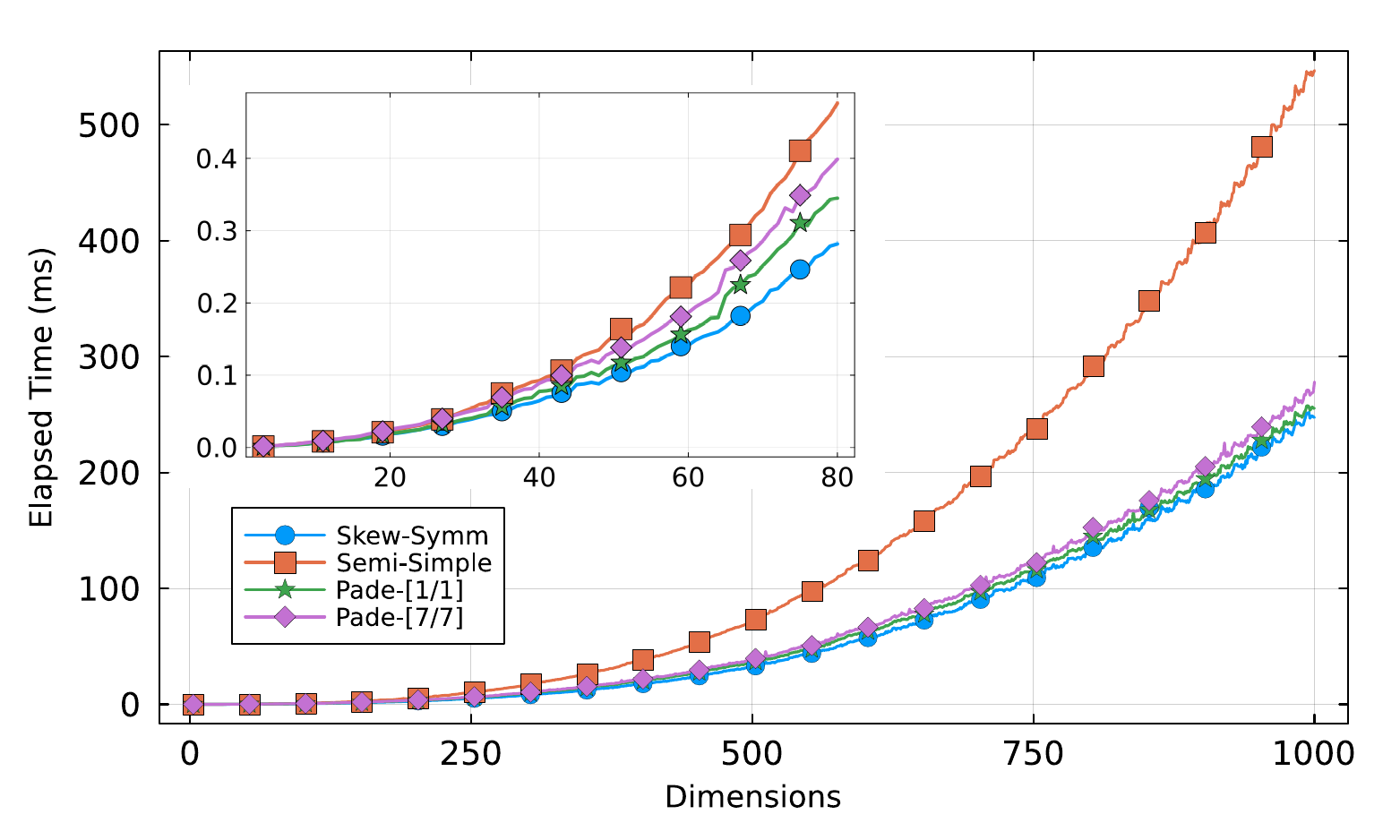}
        \caption{Total time}
        \label{fig:total_time_dlog_full}
    \end{subfigure}
    \hfill
    \begin{subfigure}[b]{0.48\textwidth}
        \centering
        \includegraphics[width=\textwidth, height=.52\textwidth]{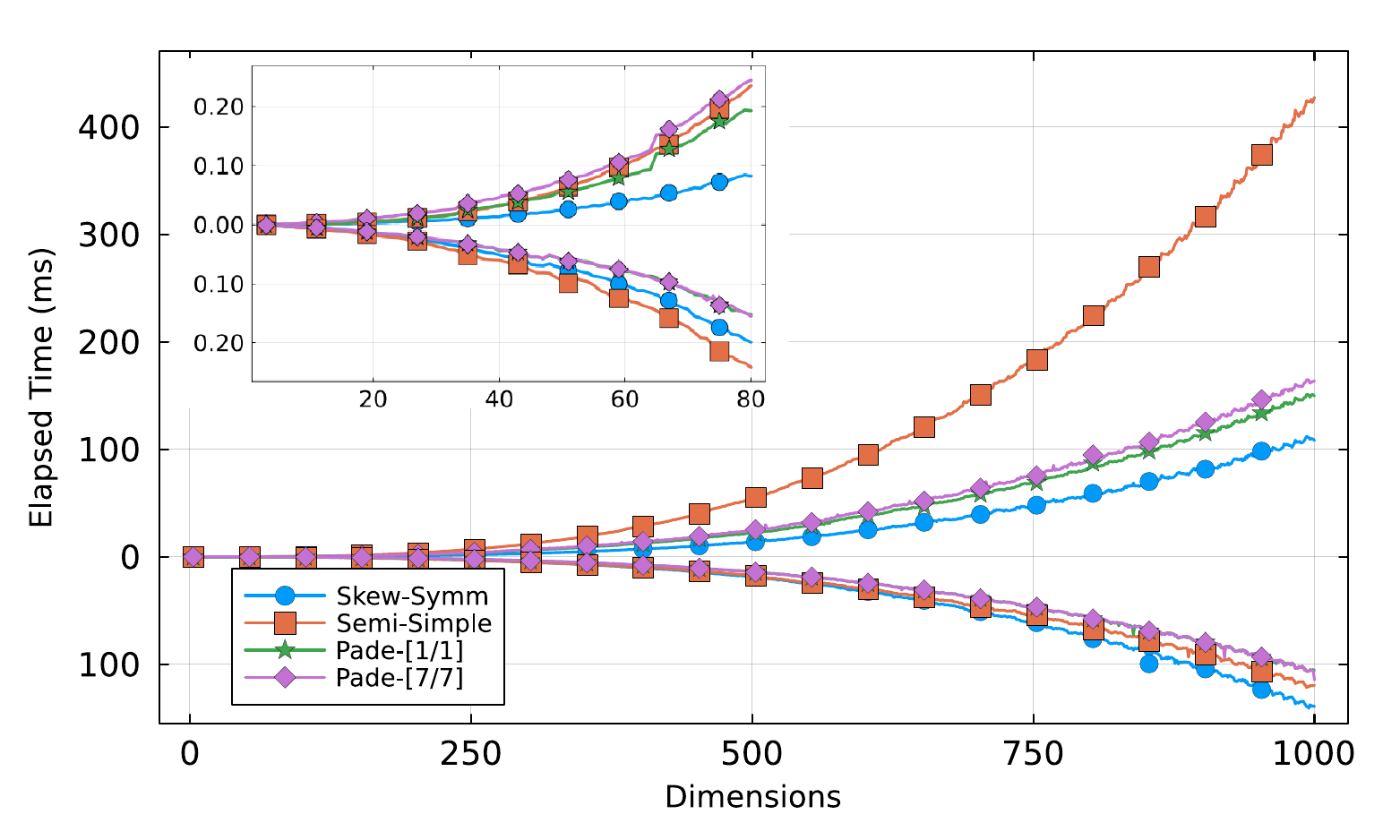}
        \caption{Breakdown (Preprocessing below $x$-axis)}
        \label{fig:breakdown_dlog_full}
    \end{subfigure}
    \caption{Execution time for $\Delta\mapsto \dexpinvof{A}[\Delta]$ across algorithms.}\label{fig:timing_dlog_full}
\end{figure}

\Cref{fig:timing_dlog_skew} and \Cref{fig:timing_dlog_full} present the efficiency comparison in $Y\mapsto X$ and $Q Y\mapsto X$ from $\dexpinvof{A}[Q Y]=X$ respectively. An excerpt of the recorded elapsed time in particular dimensions is given in \cref{tab:dlog-elapsed}.

Consider $n > 250$. For one execution of $Y\mapsto \dskew_A^{-1}(Y)$, reported in \cref{sub@fig:total_time_dlog_skew}, the proposed formula is $2.2$ times faster than Daletski$\breve{\mathrm{\i}}$--Kre$\breve{\mathrm{\i}}$n formula, and $1.3$ times faster than the less robust Pad{\'e} approximant. For the dominant computing stage in multiple calls of $Y\mapsto \dskew_A^{-1}(Y)$, reported in \cref{sub@fig:breakdown_dlog_skew}, the proposed formula is $3.6$ times faster than Daletski$\breve{\mathrm{\i}}$--Kre$\breve{\mathrm{\i}}$n formula, and $2.6$ times faster than the Pad{\'e} approximant. For one execution of $\Delta\mapsto \dexpinvof{A}[\Delta]$, reported in \cref{sub@fig:total_time_dexp_full}, the proposed formula is $3.2$ times faster than Daletski$\breve{\mathrm{\i}}$--Kre$\breve{\mathrm{\i}}$n formula, and faster than the Pad{\'e} approximant. For the dominant computing stage in multiple calls of $\Delta\mapsto \dexpinvof{A}[\Delta]$, reported in \cref{sub@fig:breakdown_dlog_full}, the proposed formula is $5.1$ times faster than Daletski$\breve{\mathrm{\i}}$--Kre$\breve{\mathrm{\i}}$n formula, and $1.7$ times faster than the order $[7/7]$ Pad{\'e} approximant.

\begin{table}[tbp]
    \centering
    \caption{Recorded Elapsed Time for $X=\dskew_A^{-1}(Y)$ and $X = \dexpinvof{A}[\Delta]$}\label{tab:dlog-elapsed}
    \tiny{
        \begin{tabular}{|l|ccc|ccc|ccc|}
            \hline
            $Y\mapsto X$& \multicolumn{3}{c|}{Preprocessing} & \multicolumn{3}{c|}{Computation} & \multicolumn{3}{c|}{Total} \\
            Problem Sizes& $10$ & $50$ & $1000$& $10$ & $50$ & $1000$& $10$ & $50$ & $1000$\\
            \hline
            Skew-Symm& 5.1e-3& 7.0e-2& \textbf{1.1e+2}& \textbf{8.2e-4}& \textbf{2.1e-2}& \textbf{9.2e+1}& \textbf{5.9e-3}& \textbf{9.1e-2}& \textbf{2.0e+2}\\
            Pad{\'e}-$[1/1]$& \textbf{5.0e-3}& \textbf{6.3e-2}& 1.3e+2& 1.4e-3& 5.8e-2& 1.8e+2& 6.4e-3& 1.2e-1& 3.1e+2\\
            Pad{\'e}-$[7/7]$& 5.0e-3& 6.4e-2& 1.4e+2& 3.7e-3& 7.7e-2& 1.9e+2& 8.7e-3& 1.4e-1& 3.3e+2\\
            Semi-Simple& 6.1e-3& 9.5e-2& 1.2e+2& 1.2e-3& 6.0e-2& 4.3e+2& 7.3e-3& 1.5e-1& 5.5e+2\\
            \hline
            $\Delta\mapsto X$& \multicolumn{3}{c|}{Preprocessing} & \multicolumn{3}{c|}{Computation} & \multicolumn{3}{c|}{Total} \\
            Problem Sizes& $10$ & $50$ & $1000$& $10$ & $50$ & $1000$& $10$ & $50$ & $1000$\\
            \hline
            Skew-Symm& 5.4e-3& 7.5e-2& 1.4e+2& \textbf{9.7e-4}& \textbf{2.4e-2}& \textbf{1.1e+2}& 6.3e-3& \textbf{9.9e-2}& \textbf{2.5e+2}\\
            Pad{\'e}-$[1/1]$& \textbf{4.8e-3}& \textbf{5.8e-2}& \textbf{1.1e+2}& 1.3e-3& 5.4e-2& 1.5e+2& \textbf{6.1e-3}& 1.1e-1& 2.6e+2\\
            Pad{\'e}-$[7/7]$& 4.8e-3& 5.9e-2& 1.1e+2& 3.6e-3& 7.4e-2& 1.6e+2& 8.4e-3& 1.3e-1& 2.8e+2\\
            Semi-Simple& 6.1e-3& 9.5e-2& 1.2e+2& 1.2e-3& 6.0e-2& 4.3e+2& 7.3e-3& 1.5e-1& 5.5e+2\\
            \hline
            \multicolumn{10}{l}{Skew-Symm: Formula~\eqref{eq:dexp-skew-inv} \:\:  Pad{\'e}-$[r/r]$: Pad{\'e} approximant~\cite{al2013logarithm} \:\: Semi-Simple: Inverse of~\eqref{eq:dexp-diag} and~\eqref{eq:linear-diag}}
        \end{tabular}
    }
\end{table}

\subsubsection{Matrix Exponential}

\begin{figure}[tbp]
    \centering
    \begin{subfigure}[b]{0.48\textwidth}
        \centering
        \includegraphics[width=\textwidth, height=.52\textwidth]{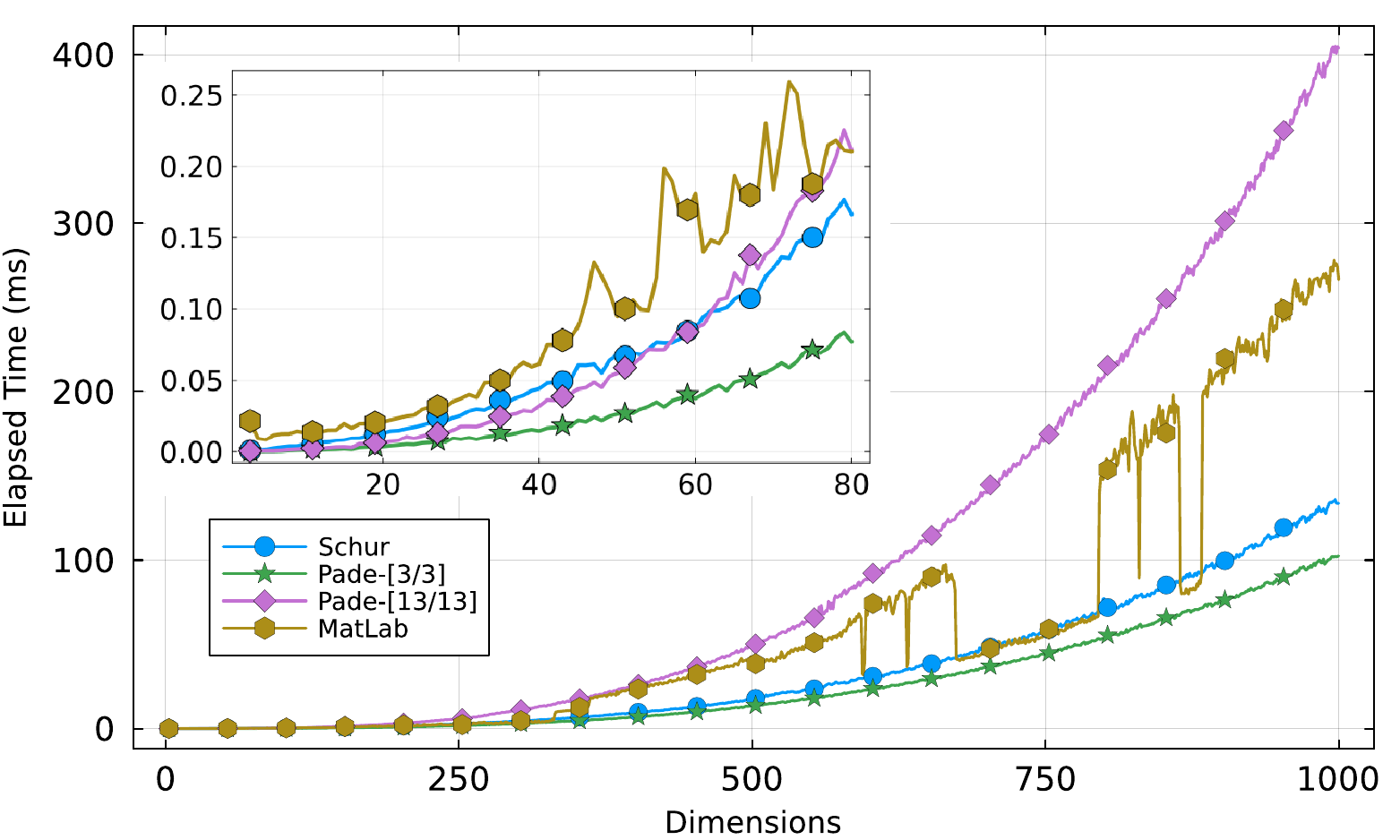}
        \caption{Compute time}
        \label{fig:expm-time}
    \end{subfigure}
    \hfill
    \begin{subfigure}[b]{0.48\textwidth}
        \centering
        \includegraphics[width=\textwidth, height=.52\textwidth]{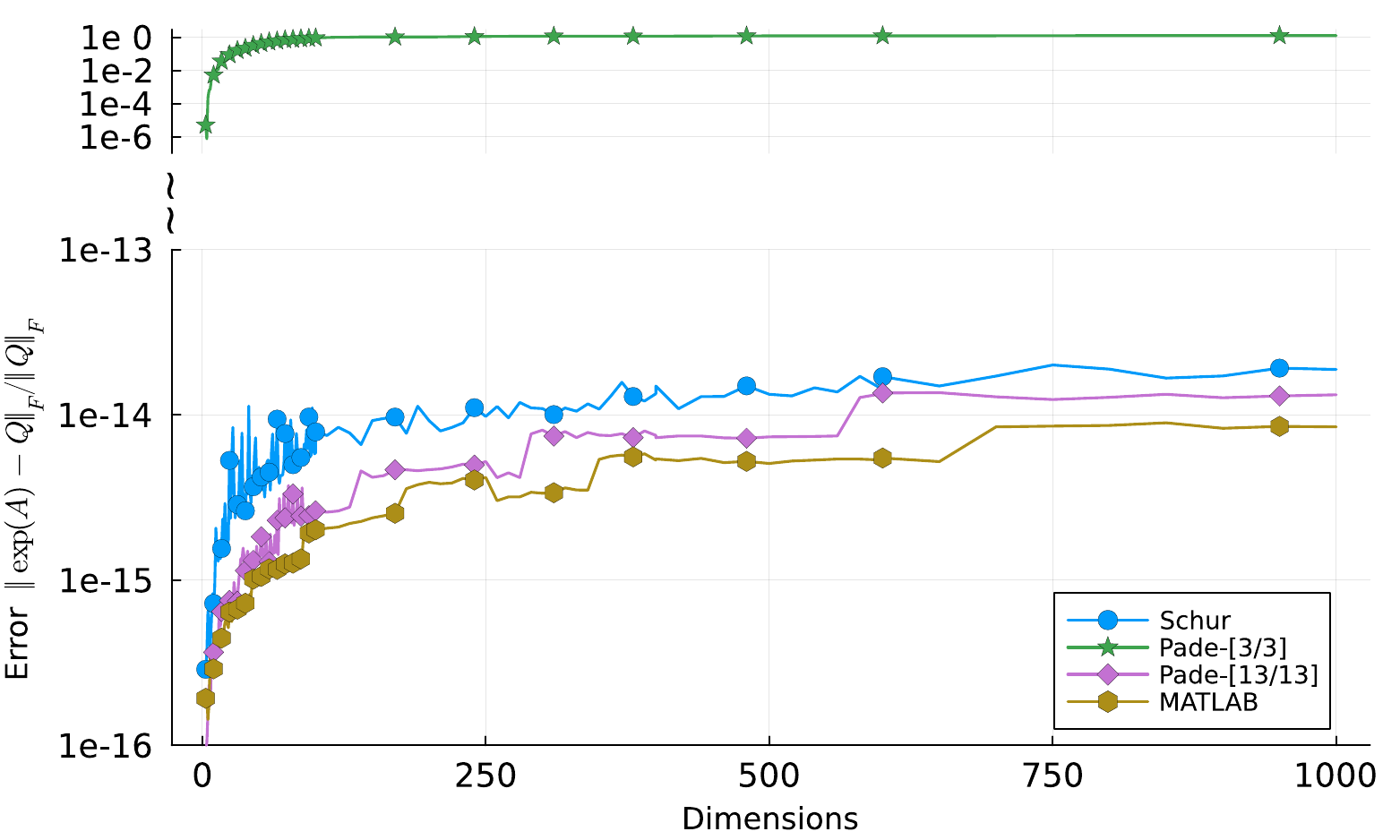}
        \caption{Compute Error}
        \label{fig:expm-error}
    \end{subfigure}
    \caption{Computation for $A\mapsto Q = \exp(A)$ across algorithms.}
    \label{fig:expm}
\end{figure}

\Cref{fig:expm} presents a comparison between different formulae for $\skewm_n\ni A\mapsto \exp(A)$. Note that the order $[3/3]$ Pad{\'e} approximant only captures at most three decimal digits of precision, while the others capture more than $16$ digits. As predicted in the complexity analysis, the exponential computed by Schur decomposition is at least faster than the \texttt{MATLAB}'s \texttt{expm}. When higher order Pad{\'e} approximant is used (by the adaptive strategy), the exponential by Schur decomposition is two times faster than \texttt{expm}, by about $100$ ms. Observe that the computation time of \texttt{expm} exhibits large jumps owing to the adaptive Pad{\'e} order. Unlike a fixed-order Pad{\'e} approximant, \texttt{MATLAB} dynamically determines the Pad{\'e} order for computing $\exp(A)$ based on the input matrix $A$, which explains these discontinuities in elapsed time. While the details are not publicly documented, the development history of the improved single-precision \texttt{expm} (R2023b) in the MathWorks documentation \cite{MATLABexpm} suggests that a dynamic Pad{\'e} order is used internally. Note that the elapsed time of \texttt{expm} lies between those of the lowest $[3/3]$ and the highest $[13/13]$ orders reported in \cite{al2009exponential}, which is consistent with typical behavior of a dynamically selected Pad{\'e} strategy.

\section{Conclusion}~\label{sec:conclusion}
This paper investigates the differentiation of the matrix exponential restricted to the skew-symmetric matrices. In particular, it derives the efficient formulae for the differentiation~\eqref{eq:dskew-skew} and its inverse~\eqref{eq:dexp-skew-inv} where computations are performed in real arithmetic. The necessary and sufficient condition for the invertible differentiation is identified in~\cref{thm:invertibility-skew-symm}, which is used to characterize the distance to rank-deficient differentiation in~\cref{thm:distance-bound-to-conjugate-locus}. These new results lead to the definition of nearby matrix logarithm as a well-defined local inversion to the matrix exponential in~\cref{thm:diffeomorphism}, which subsumes the principal logarithm and works beyond the image of the principal logarithm. The preliminary experiments given in \Cref{sec:nearby-logarithm} have demonstrated the potential of the nearby logarithm to recognize continuous structures. Finally, the derived symbolic formulae have been efficiently implemented in \texttt{C++} and demonstrated significant computational advantages over existing formulae.

\bibliographystyle{alphaurl}
\bibliography{refs}

@article{edelman1998,
  author = {Edelman, Alan and Arias, Tom\'{a}s A. and Smith, Steven T.},
  title = {{T}he {G}eometry of {A}lgorithms with {O}rthogonality {C}onstraints},
  journal = {SIAM Journal on Matrix Analysis and Applications},
  volume = {20},
  number = {2},
  pages = {303-353},
  year = {1998},
  publisher = {SIAM}
}

@article{jurdjevic2020extremal,
  title     = {{E}xtremal {C}urves on {S}tiefel and {G}rassmann {M}anifolds},
  author    = {Jurdjevic, Velimir and Markina, Irina and Leite, F Silva},
  journal   = {The Journal of Geometric Analysis},
  volume    = {30},
  number    = {4},
  pages     = {3948--3978},
  year      = {2020},
  publisher = {Springer}
}

@article{dieci1999real,
  title     = {{O}n {R}eal {L}ogarithms of {N}earby {M}atrices and {S}tructured {M}atrix {I}nterpolation},
  author    = {Dieci, Luca and Morini, Benedetta and Papini, Alessandra and Pasquali, Aldo},
  journal   = {Applied Numerical Mathematics},
  volume    = {29},
  number    = {1},
  pages     = {145--165},
  year      = {1999},
  publisher = {Elsevier}
}

@article{najfeld1995,
  title     = {{D}erivatives of the {M}atrix {E}xponential and {T}heir {C}omputation},
  author    = {Najfeld, Igor and Havel, Timothy F},
  journal   = {Advances in Applied Mathematics},
  volume    = {16},
  number    = {3},
  pages     = {321},
  year      = {1995},
  publisher = {New York, Academic Press.}
}

@article{zimmermann2021computing,
  author = {Zimmermann, Ralf and H\"{u}per, Knut},
  title = {{C}omputing the {R}iemannian {L}ogarithm on the {S}tiefel {M}anifold: {M}etrics, {M}ethods, and {P}erformance},
  journal = {SIAM Journal on Matrix Analysis and Applications},
  volume = {43},
  number = {2},
  pages = {953-980},
  year = {2022},
  publisher = {SIAM}
}

@article{al2013logarithm,
  title={{C}omputing the {F}r{\'e}chet {D}erivative of the {M}atrix {L}ogarithm and {E}stimating the {C}ondition {N}umber},
  author={Al-Mohy, Awad H and Higham, Nicholas J and Relton, Samuel D},
  journal={SIAM Journal on Scientific Computing},
  volume={35},
  number={4},
  pages={C394--C410},
  year={2013},
  publisher={SIAM}
}

@article{al2009exponential,
  title={{C}omputing the {F}r{\'e}chet {D}erivative of the {M}atrix {E}xponential, {W}ith an {A}pplication to {C}ondition {N}umber {E}stimation},
  author={Al-Mohy, Awad H and Higham, Nicholas J},
  journal={SIAM Journal on Matrix Analysis and Applications},
  volume={30},
  number={4},
  pages={1639--1657},
  year={2009},
  publisher={SIAM}
}

@book{golub2013matrix,
  title={{M}atrix {C}omputations},
  author={Golub, Gene H and Van Loan, Charles F},
  year={2013},
  publisher={JHU press}
}

@book{rossmann2006lie,
  title     = {{L}ie {G}roups: {A}n {I}ntroduction {T}hrough {L}inear {G}roups},
  author    = {Rossmann, Wulf},
  volume    = {5},
  year      = {2006},
  publisher = {Oxford University Press on Demand}
}

@article{higham2005scaling,
  author = {Higham, Nicholas J.},
  title = {{T}he {S}caling and {S}quaring {M}ethod for the {M}atrix {E}xponential {R}evisited},
  journal = {SIAM Journal on Matrix Analysis and Applications},
  volume = {26},
  number = {4},
  pages = {1179-1193},
  year = {2005},
  publisher = {SIAM}
}

@article{al2010new,
  title={{A} {N}ew {S}caling and {S}quaring {A}lgorithm for the {M}atrix {E}xponential},
  author={Al-Mohy, Awad H and Higham, Nicholas J},
  journal={SIAM Journal on Matrix Analysis and Applications},
  volume={31},
  number={3},
  pages={970--989},
  year={2010},
  publisher={SIAM}
}

@article{turaga2011statistical,
  title={{S}tatistical {C}omputations on {G}rassmann and {S}tiefel {M}anifolds for {I}mage and {V}ideo-{B}ased {R}ecognition},
  author={Turaga, Pavan and Veeraraghavan, Ashok and Srivastava, Anuj and Chellappa, Rama},
  journal={IEEE Transactions on Pattern Analysis and Machine Intelligence},
  volume={33},
  number={11},
  pages={2273--2286},
  year={2011},
  publisher={IEEE}
}

@article{wilcox1967exponential,
  title={{E}xponential {O}perators and {P}arameter {D}ifferentiation in {Q}uantum {P}hysics},
  author={Wilcox, Ralph M},
  journal={Journal of Mathematical Physics},
  volume={8},
  number={4},
  pages={962--982},
  year={1967},
  publisher={American Institute of Physics}
}

@article{chan2003generalized,
  title={{A} {G}eneralized {L}inear {M}odel for {R}epeated {O}rdered {C}ategorical {R}esponse {D}ata},
  author={Chan, KS and Munoz-Hernandez, B},
  journal={Statistica Sinica},
  pages={207--226},
  year={2003},
  publisher={JSTOR}
}

@article{chen2001analytic,
  title={{A}nalytic {D}erivatives of the {M}atrix {E}xponential for {E}stimation of {L}inear {C}ontinuous-{T}ime {M}odels},
  author={Chen, Baoline and Zadrozny, Peter A},
  journal={Journal of Economic Dynamics and Control},
  volume={25},
  number={12},
  pages={1867--1879},
  year={2001},
  publisher={Elsevier}
}

@article{mataigne2024eigenvalue,
  title={{T}he {E}igenvalue {D}ecomposition of {N}ormal {M}atrices by the {D}ecomposition of the {S}kew-{S}ymmetric {P}art {W}ith {A}pplications to {O}rthogonal {M}atrices},
  volume={42}, 
  url={https://journals.uwyo.edu/index.php/ela/article/view/9957}, 
  DOI={10.13001/ela.2026.9957}, 
  number={42}, 
  journal={The Electronic Journal of Linear Algebra}, 
  author={Mataigne, Simon and Gallivan, Kyle A.}, 
  year={2026}, 
  month={May}, 
  pages={349–386}
}

@inproceedings{zhou2019continuity,
  title={{O}n the {C}ontinuity of {R}otation {R}epresentations in {N}eural {N}etworks},
  author={Zhou, Yi and Barnes, Connelly and Lu, Jingwan and Yang, Jimei and Li, Hao},
  booktitle={Proceedings of the IEEE/CVF Conference on Computer Vision and Pattern Recognition},
  pages={5745--5753},
  year={2019}
}

@inproceedings{wang2020orthogonal,
  title={{O}rthogonal {C}onvolutional {N}eural {N}etworks},
  author={Wang, Jiayun and Chen, Yubei and Chakraborty, Rudrasis and Yu, Stella X},
  booktitle={Proceedings of the IEEE/CVF Conference on Computer Vision and Pattern Recognition},
  pages={11505--11515},
  year={2020}
}

@article{qi2021orthogonalmoment,
  title={{A} {S}urvey of {O}rthogonal {M}oments for {I}mage {R}epresentation: {T}heory, {I}mplementation, and {E}valuation},
  author={Qi, Shuren and Zhang, Yushu and Wang, Chao and Zhou, Jiantao and Cao, Xiaochun},
  journal={ACM Computing Surveys (CSUR)},
  volume={55},
  number={1},
  pages={1--35},
  year={2021},
  publisher={ACM New York, NY}
}

@article{gallier2003computing,
  title={{C}omputing {E}xponentials of {S}kew-{S}ymmetric {M}atrices and {L}ogarithms of {O}rthogonal {M}atrices},
  author={Gallier, Jean and Xu, Dianna},
  journal={International Journal of Robotics and Automation},
  volume={18},
  number={1},
  pages={10--20},
  year={2003},
  publisher={Anaheim, Calif.; Calgary, Alta.: Acta Press,[1986]-}
}

@book{brigham1988fast,
author = {Brigham, E. Oran},
title = {{T}he {F}ast {F}ourier {T}ransform and {I}ts {A}pplications},
year = {1988},
isbn = {0133075052},
publisher = {Prentice-Hall, Inc.},
}

@article{daletskii1965integration,
  title={{I}ntegration and {D}ifferentiation of {F}unctions of {H}ermitian {O}perators and {A}pplications to the {T}heory of {P}erturbations},
  author={Daletski{\v{i}}, Ju L and Kre{\v{i}}n, Selim Grigorievich},
  journal={AMS Translations (2)},
  volume={47},
  number={1-30},
  pages={10--1090},
  year={1965}
}

@article{meyer1986continuous,
  title={{C}ontinuous {O}rthonormalization for {B}oundary {V}alue {P}roblems},
  author={Meyer, Gunter H},
  journal={Journal of Computational Physics},
  volume={62},
  number={1},
  pages={248--262},
  year={1986},
  publisher={Elsevier}
}

@inproceedings{dieci2000continuous,
  title={{C}ontinuous {O}rthonormalization for {L}inear {T}wo-{P}oint {B}oundary {V}alue {P}roblems {R}evisited},
  author={Dieci, Luca and van Vleck, Erik S},
  booktitle={Dynamics of Algorithms},
  pages={69--90},
  year={2000},
  organization={Springer}
}

@article{dieci1997compuation,
  title={{O}n the {C}ompuation of {L}yapunov {E}xponents for {C}ontinuous {D}ynamical {S}ystems},
  author={Dieci, Luca and Russell, Robert D and Van Vleck, Erik S},
  journal={SIAM Journal on Numerical Analysis},
  volume={34},
  number={1},
  pages={402--423},
  year={1997},
  publisher={SIAM}
}

@article{barker2018evans,
  title={{E}vans {F}unction {C}omputation for the {S}tability of {T}ravelling {W}aves},
  author={Barker, Blake and Humpherys, Jeffrey and Lyng, Gregory and Lytle, Joshua},
  journal={Philosophical Transactions of the Royal Society A: Mathematical, Physical and Engineering Sciences},
  volume={376},
  number={2117},
  pages={20170184},
  year={2018},
  publisher={The Royal Society Publishing}
}

@article{borzsak1996lyapunov,
  title = {{L}yapunov {I}nstability of {F}luids {C}omposed of {R}igid {D}iatomic {M}olecules},
  author = {Borzs\'ak, Istv\'an and Posch, H. A. and Baranyai, Andr\'as},
  journal = {Phys. Rev. E},
  volume = {53},
  issue = {4},
  pages = {3694--3701},
  numpages = {0},
  year = {1996},
  month = {Apr},
  publisher = {American Physical Society}
}

@book{horn2012matrix,
  title={{M}atrix {A}nalysis},
  author={Horn, Roger A and Johnson, Charles R},
  year={2012},
  publisher={Cambridge University Press}
}

@phdthesis{deng2024smoothly,
  title={{S}moothly {E}volving {G}eodesics in the {S}pecial {O}rthogonal {G}roup: {D}efinitions, {C}omputations and {A}pplications},
  author={Deng, Zhifeng},
  year={2024},
  school={The Florida State University}
}

@article{absil2025ultimate,
author = {Absil, P.-A. and Mataigne, Simon},
title = {{T}he {U}ltimate {U}pper {B}ound on the {I}njectivity {R}adius of the {S}tiefel {M}anifold},
journal = {SIAM Journal on Matrix Analysis and Applications},
volume = {46},
number = {2},
pages = {1145-1167},
year = {2025},
}

@article{cardoso2010exponentials,
  title={{E}xponentials of {S}kew-{S}ymmetric {M}atrices and {L}ogarithms of {O}rthogonal {M}atrices},
  author={Cardoso, Jo{\~a}o R and Leite, F Silva},
  journal={Journal of Computational and Applied Mathematics},
  volume={233},
  number={11},
  pages={2867--2875},
  year={2010},
  publisher={Elsevier}
}

@article{GuiguiMiolanePennec2023,
url_hide = {http://dx.doi.org/10.1561/2200000098},
year = {2023},
volume = {16},
journal = {Foundations and Trends in Machine Learning},
title = {{I}ntroduction to {Riemannian} {G}eometry and {G}eometric {S}tatistics: {F}rom {B}asic {T}heory to {I}mplementation with {Geomstats}},
doi = {10.1561/2200000098},
issn = {1935-8237},
number = {3},
pages = {329-493},
author = {Nicolas Guigui and Nina Miolane and Xavier Pennec},
}

@book {GallierQuaintance2020,
    AUTHOR = {Gallier, Jean and Quaintance, Jocelyn},
     TITLE = {{D}ifferential {G}eometry and {L}ie {G}roups---{A} {C}omputational
              {P}erspective},
    SERIES = {Geometry and Computing},
    VOLUME = {12},
 PUBLISHER = {Springer, Cham},
      YEAR_orig = {[2020] \copyright 2020},
      year = 2020,
     PAGES = {xv+777},
      ISBN = {978-3-030-46039-6; 978-3-030-46040-2},
   MRCLASS = {53-01 (53C30 53C35)},
  MRNUMBER = {4164718},
  doi = "10.1007/978-3-030-46040-2",
}

@misc{MATLABexpm,
  author       = {MathWorks},
  title        = {expm: Matrix Exponential},
  howpublished = {\url{https://www.mathworks.com/help/matlab/ref/expm.html}},
  note         = {Accessed: 2025-11-06},
  year         = {2025}
}
\end{document}